\documentclass[11pt]{article}
\usepackage{amssymb,latexsym}
\usepackage{epsfig}
\usepackage{eufrak}
\usepackage{amsmath}
\usepackage{mathrsfs}
\usepackage{color}

\newtheorem{theorem}{Theorem}
\newtheorem{lemma}{Lemma}
\newtheorem{corollary}{Corollary}

\newcommand{\be}{\begin{equation}}
\newcommand{\ee}{\end{equation}}
\newcommand{\bee}{\begin{eqnarray*}}
\newcommand{\eee}{\end{eqnarray*}}
\newcommand{\bel}{\begin{eqnarray}}
\newcommand{\eel}{\end{eqnarray}}
\newcommand{\bec}{\begin{cases}}
\newcommand{\eec}{\end{cases}}
\newcommand{\bem}{\begin{bmatrix}}
\newcommand{\eem}{\end{bmatrix}}

\newcommand{\la}{\label}
\newcommand{\li}{\left}
\newcommand{\ri}{\right}

\newcommand{\ovl}{\overline}
\newcommand{\udl}{\underline}

\newcommand{\vep}{\varepsilon}
\newcommand{\lm}{\lambda}

\newcommand{\si}{\sigma}

\newcommand{\de}{\delta}

\newcommand{\vse}{\vartheta}
\newcommand{\se}{\theta}

\newcommand{\ze}{\zeta}
\newcommand{\al}{\alpha}

\newcommand{\f}{\frac}
\newcommand{\sq}{\sqrt}
\newcommand{\cd}{\cdots}

\newcommand{\mscr}{\mathscr}

\newcommand{\bb}{\mathbb}

\newcommand{\wh}{\widehat}

\newcommand{\sh}{\slash}

\newcommand{\tx}{\text}

\newcommand{\iy}{\infty}

\newcommand{\bed}{\begin{description}}
\newcommand{\eed}{\end{description}}
\newcommand{\bei}{\begin{itemize}}
\newcommand{\eei}{\end{itemize}}
\newcommand{\ben}{\begin{enumerate}}
\newcommand{\een}{\end{enumerate}}
\newcommand{\bib}{\bibitem}
\newcommand{\beL}{\begin{lemma}}
\newcommand{\eeL}{\end{lemma}}
\newcommand{\beT}{\begin{theorem}}
\newcommand{\eeT}{\end{theorem}}

\newcommand{\bpf}{\begin{pf}}
\newcommand{\epf}{\end{pf}}

\newcommand{\pfbox}{\hfill\mbox{$\Box$}}

\newenvironment{pf}{\paragraph*{Proof{\rm.}}}{\pfbox\bigskip}

\begin{document}

\title{{\bf A Statistical Theory for Measurement and Estimation
\\ of Rayleigh Fading Channels}
\thanks{The authors are with
Department of Electrical and Computer Engineering,
Louisiana State University, Baton Rouge, LA 70803; Email:
\{chan,ggu, kemin\}@ece.lsu.edu, Tel: (225)578-\{8961, 5534,5533\}, and
Fax: (225) 578-5200. }}

\author{Xinjia Chen, \ Guoxiang Gu, \ and \ Kemin Zhou}

\date{June 2007}

\maketitle

\begin{abstract}

In this paper, we propose a statistical theory on measurement and
estimation of Rayleigh fading channels in wireless communications
and provide complete solutions to the fundamental problems: What is
the optimum estimator for the statistical parameters associated with
the Rayleigh fading channel, and how many measurements are
sufficient to estimate these parameters with the prescribed margin
of error and confidence level? Our proposed statistical theory
suggests that two testing signals of different strength be used. The
maximum likelihood (ML) estimator is obtained for estimation of the
statistical parameters of the Rayleigh fading channel that is both
sufficient and complete statistic. Moreover, the ML estimator is the
minimum variance (MV) estimator that in fact achieves the
Cram\'{e}r-Rao lower bound.

\end{abstract}

\section{Introduction}

In mobile radio channels, the Rayleigh distribution is commonly used
to describe the statistical nature of the received envelope of a
flat fading signal, or the envelope of an individual multipath
component. Flat fading is often associated with the narrow band
channel. By assuming that the real and imaginary parts of the
channel gain are independent Gaussian random variables with zero
mean and equal variance, the amplitude of the channel gain or PL
(path loss) becomes a Rayleigh random variable. For wide band
channels, multipath gains are typically assumed to be uncorrelated
scattering (US) and each path gain is wide-sense stationary (WSS)
which are termed as WSSUS channels \cite{b63}. By assuming that the
real and imaginary parts of each path gain are independently
Gaussian distributed with zero mean and equal variance, the
amplitudes of the multipath gains become independent Rayleigh random
variables. Theoretically a Rayleigh random variable is uniquely
specified by its 2nd moment that is the sum of the variances of its
two independent Gaussian components. In practice the statistics of
the Rayleigh fading channel are incomplete without additional
knowledge of the noise power and SNR (signal-to-noise ratio). That
is, the statistical parameters of the channel power, noise power,
and SNR together characterize the Rayleigh fading channel
completely, giving rise to the measurement and estimation problem
for these statistical parameters.

\vspace{1mm}

Estimation for the 2nd moment of a Rayleigh random variable has
practical importance in channel modeling and estimation \cite{besw,
czink, gghn}, and in radio coverage, location, and measurement
\cite{b, c}. For these reasons statistical estimation of Rayleigh
fading channels has been studied and reported in the existing
research literature. The early work of Peritsky \cite{per} shows
that the simple averaged square of the signal strength based on
i.i.d. (independent and identically distributed) measurement samples
is both an ML (maximum likelihood) and MV (minimum variance)
estimator, and such an averaged square is a sufficient and complete
statistic. However his results focus only on the noise-free case and
have limited applications. The same problem has been investigated in
\cite{ak, cb, chen, dietrich, dogandzic, tag, valenzuela} that
encompass the Ricean and Nakagami fading distributions, as well as
MIMO fading channels. But there lack optimum estimators that achieve
the Cram\'{e}r-Rao lower bound.

\vspace{1mm}

In this
paper we propose a statistical theory
on measurement and estimation of Rayleigh fading channels.
Our contributions include derivation of explicit
{\it a priori} bounds on the measurement sample size
that achieve the
prescribed margin of error and
the confidence level, and discovery of the optimum estimator
that is both an ML and MV estimator, and
achieves the Cram\'{e}r-Rao lower bound.
These results complement the existing work reported in the literature,
In the noise-free case our sample
size bound resembles the one implicitly
indicated in \cite{chernoff}
which asserts that,
to estimate the binomial probability with the
prescribed margin of absolute error $\vep$ and confidence level $1 -
\de$, it suffices to have a sample size greater than {\tiny $\f{\ln
(\f{2}{\de})}{ 2 \vep^2 }$}. One notable difference is that our
sample size bound is derived for the relative error while the bound
in \cite{chernoff} is derived for the absolute error in estimation
of the binomial probability.
Furthermore an interval estimate similar to that
in \cite{per} is obtained with much simpler calculation.

\vspace{1mm}

In applications to Rayleigh fading channels, noisy measurement
samples have to be taken into consideration. Our proposed
statistical theory suggests that two testing signals of different
strength be used. The use of two different testing signals enables
us to extend the sample complexity results from the noiseless case
to the noisy case which solves the sample complexity problem not
only for the 2nd moment or the mean channel power, but also for the
mean noise power and for the SNR. The explicit {\it a priori} bounds
resemble to those in the noise-free case of which the sample size is
roughly inversely proportional to the square of the margin of error
and is linear with respect to the logarithm of the inverse of the
gap between the confidence level and $1$. Our results show that a
typical margin of error can be achieved with near certainty and
modest sample size. More importantly an optimum estimator is
obtained for noisy Rayleigh fading channels that is both an ML and
MV estimator. It inherits the same sufficient and complete statistic
property from that in the noiseless case \cite{per}, and it achieves
the Cram\'{e}r-Rao lower bound.

\vspace{1mm}

The content of the paper is organized as follows. In Section 2 we
first present the Rayleigh fading channel and its associated
measurement and estimation problem in terms of the sample complexity
and the optimum estimator. The sample
complexity problem is then addressed in Section 3 by establishing
an explicit {\it a priori} bound on the sample size
that is asymptotically tight based
on noiseless measurements. In Section 4 we present our
main results on parameter estimation for Rayleigh fading channels
by
taking the noisy measurement samples into consideration
to which complete solutions are derived for measurement and estimation
of Rayleigh fading channels.
Numerical examples are presented in Section 5 to
illustrate our results, and the paper is concluded in
Section 6. The notations are
standard and will be made clear as we proceed.

\section{Rayleigh Fading Channels}

In a typical urban environment, there is no LOS (line of sight)
between the transmitter and the receiver. The wireless channel is
characterized by multipath and follows the Rayleigh distribution.
Specifically a widely used channel model is the following
continuous-time CIR (channel impulse response): \be c(t; \tau) =
\sum_{n=0}^N
\alpha_n(t)e^{j\beta_n(t)}\delta\left(t-\tau_n(t)\right) \ee where
$\alpha_n(t)$, $\beta_n(t)$, and $\tau_n(t)$ are the amplitude,
phase, and time delay, respectively, associated with the $n$th path,
and $\delta(\cdot)$ is the Dirac delta function. For narrow band
systems, the channel gain is given by \be H(t) = H_R(t)+jH_I(t):=
\left[\sum_{n=1}^N \alpha_n(t)\cos(\phi_n(t))\right]
+j\left[\sum_{n=1}^N \alpha_n(t)\sin(\phi_n(t))\right] \ee The
arguments $\{\phi_n(t)\}$ are due to $\{\beta_n(t)\}$ and
$\{\tau_n(t)\}$ induced by Doppler shifts and time delays that form
a set of independent random variables uniformly distributed over
$[-\pi, \; \pi]$. Under the assumption that $\alpha_n(t)\equiv
\alpha_n$, $\beta_n(t)\equiv \beta_n$, and $\tau_n(t)\equiv \tau_n$,
$H(t)$ becomes a WSS (wide-sense stationary) process. For large $N$,
the central limiting theorem can be invoked to treat both $H_R(t)$
and $H_I(t)$ as independent Gaussian random processes with zero mean
and equal variance \cite{p, st}. This gives rise to the Rayleigh
fading channel in light of the fact that $|H(t)|^2$ is a Rayleigh
random variable. For wide band systems, multipath gains can be
resolved up to the resolution dictated by the sampling frequency
$f_s$. An equivalent discretized CIR can be represented by \be H(t)
= \sum_{k=0}^{L-1} H_k(t)\delta_{t-kT_s}, \; \; \; \; \;
T_s=f_s^{-1} \ee in which $|H_k(t)|^2$ is Rayleigh distributed for
$0\leq k< L$ where $\delta_t$ is the Kronekar delta function.

\vspace{1mm}

Rayleigh fading channels are widely used in wireless communications.
An important and practical problem is measurement
of $|H(t)|^2$ and estimation of $\bb{E}[|H(t)|^2]$,
with $\bb{E}[\ \cdot\ ]$ the expectation, for narrow band systems,
or $|H_k(t)|^2$ and power delay profile $\bb{E}[|H_k(t)|^2]$ for
wide band systems. Such a problem is crucial in
modeling and estimation of wireless channels, and
in radio coverage and location. A common strategy is
to transmit i.i.d. pilot symbols or testing signals
to obtain measurement samples of $|H(t)|^2$
or $|H_k(t)|^2$, and then to estimate the mean. For simplicity
we drop the time argument $t$, and focus on
narrow band systems but our results
are applicable to wide band systems as well by converting
frequency-selective fading channels into flat fading channels
using the OFDM (orthogonal frequency division multiplexing)
scheme \cite{gghn}.
The detail is omitted.

\vspace{1mm}

Denote $\Re[C]$ as the real part of the complex number $C$. Let the
transmitted symbol be $S$. In the complex form, the received signal
can be written as (``$\Longrightarrow$'' stands for ``implying'')
\be x(t) = \Re [ (HS + V)e^{j 2\pi f_ct}] \; \; \; \;
\Longrightarrow \; \; \; \; X = HS + V \ee where $f_c$ is the
carrier frequency and $V = V_R + jV_I$ is the complex additive
Gaussian noise with zero mean and equal variance $\sigma_V^2$. Since
the average power of the transmitted signal is $P_s = \bb{E}[|S|^2]$
and the average power of the received signal is $\bb{E}[|HS|^2]$,
the channel power or PL is
\[
{\rm PL} = 10 \log_{10}\left(\f{ \bb{E}[|S|^2] } { \bb{E}[|HS|^2]}\right) =
10
\log_{10}\left(\f{ \bb{E}[|S|^2] } { \bb{E}[|H|^2] \; \bb{E} [|S|^2]}\right)
= -
10 \log_{10}\left(\bb{E} [|H|^2]\right) = - 10 \log_{10} \left(2
\si_H^2\right)
\]
where $\si_H^2$ is the common variance of $H_R$ and $H_I$. It can be
seen that $\si_H^2$ depends on the attenuation factors $\{\al_n\}$
that decrease as the receiver is further away from the transmitter.
It follows that $\si_H^2$ and thus SNR become small when the
receiver is far away from the transmitter. To predict the coverage
and the performance of wireless systems, it is important to know the
operating range of the SNR at different locations. The SNR at the
receiver is easily found to be
\[
{\rm SNR} = 10 \log_{10} \left(\f{ \bb{E}[|HS|^2]} { \bb{E}[|V|^2] }\right)
= 10
\log_{10}\left(\f{ \bb{E}[|S|^2] \; \bb{E}[|H|^2]} { \bb{E}[|V|^2] }\right)
= 10
\log_{10}\left(\bb{E}[|S|^2]\right) + 10 \log_{10} \li ( \f{ \si_H^2 } {
\si_V^2
} \ri )
\]
where $\si_V^2$ is the common variance of $V_R$ and $V_I$.
Since
$P_s =
\bb{E}[|S|^2]$, the average power for transmitting
the symbol $S$, is computable, to obtain an estimate
for the SNR it is necessary to
have an estimate for the ratio of $\si_H^2$ to $\si_V^2$.
However it remains unclear whether or not the knowledge of SNR
helps to estimate the PL that will be investigated later.
Our main results to be presented
show that with appropriate signaling, the variances
$\si_H^2$, $\si_V^2$, and the SNR can all be estimated
in satisfaction.

\vspace{1mm}

The objective of this paper is in fact beyond estimation of the PL
that has been investigated by several researchers \cite{chen,
dietrich, dogandzic, per, tag, valenzuela}. Our goal is to develop a
relevant statistical theory in order to derive an optimum estimator
for the statistical parameters associated with the Rayleigh fading
channel, and to derive {\it a priori} sample complexity bounds given
the prescribed margin of error and the confidence level. These two
problems are not completely solved in the existing research
literature that will be studied thoroughly in the next several
sections. We will begin with the ideal case of noiseless
measurements by assuming independent PSK (phase shift keying) or
BPSK (binary PSK) signaling at the transmitter leading to $N$
independent measurement samples $\{X_i\}_{i=1}^N$ at the receiver.
In light of \cite{per}, \be \widehat{\sigma_H^2} =
\frac{\overline{X_N^2}}{P_s}, \; \; \; \; \; \overline{X_N^2} =
\frac{1}{2N}\sum_{i=1}^N |X_i|^2 \ee is an ML estimator for
$\sigma^2_H$ that is both statistically sufficient and complete. We
thus have an equivalent problem for estimation of the 2nd moment of
the Rayleigh random variable based on its $N$ i.i.d. samples. The
complete results on sample complexity in this case will be presented
in the next section.

\vspace{1mm}

After solving the sample complexity problem in the noise-free
case, we will investigate the sample complexity problem
in estimation of the PL, the noise power, and the SNR,
associated with Rayleigh fading channels
based on noisy measurement samples. Two sets of
independent BPSK signals $\{S_{1,i}\}_{i=1}^{N_1}$
and $\{S_{2,i}\}_{i=1}^{N_2}$ are transmitted with
difference in their average power: $|S_{1,i}|^2=P_{s1}$
and $|S_{2,i}|^2=P_{s2}$ for all possible $i$ where
$P_{s1}>P_{s2}\geq 0$. The following two quantities
\be\label{tmp0}
\overline{X_{1,N_1}^2} = \frac{1}{2N_1}\sum_{i=1}^{N_1}
|X_{1,i}|^2, \; \; \; \; \;
\overline{X_{2,N_2}^2} = \frac{1}{2N_2}\sum_{i=1}^{N_2}
|X_{2,i}|^2
\ee
turn out to be sufficient and complete statistics
in statistical estimation of the
Rayleigh fading channel based on noisy samples.
In fact the optimum estimator is obtained
that achieves the Cram\'{e}r-Rao lower bound.
Complete results are obtained and will be presented in a later section.

\section{Sample Complexity in the Noiseless Case}

As discussed earlier estimation of the 2nd moment for
the Rayleigh fading channel in the noise-free case is equivalent
to that for the Rayleigh random variable.
Let $X$ be a Rayleigh random variable with PDF
(probability density function)
\be
f_X(x) = \bec \f{2 x}{\mu} e^{- \f{x^2}{\mu}}, & \; \;
\tx{if} \; \; \; x \geq 0\\
\; \; \; \; \; 0, &  \; \; \tx{if} \; \; \; x < 0 \eec \ee Then
$X^2={X_R^2+X_I^2}$ where $X_R$ and $X_I$ are independent Gaussian
random variables with zero mean and equal variance $\sigma^2=\mu/2$.
As it is seen, the PDF of the Rayleigh random variable is completely
specified by its 2nd moment $\mu=\bb{E}[X^2] = 2\sigma^2$. To
estimate $\mu$, it is conventional to obtain $N$ i.i.d. observations
$X_1, \cd, X_N$ and take
\be \label{muhat} \wh{\mu}_N = \f{X_1^2 +
X_2^2 + \cdots + X_N^2}{N} \ee
as the estimator of $\mu$. It is
shown in \cite{per} that the above estimate is a sufficient
statistic with good performance and asymptotically unbiased with
large sample size. In addition an algorithm is derived in \cite{per}
to compute the interval estimate of $\mu$ that is a function of the
sample size and involves the use of the tabulated $\chi^2$ CDF.
However thus far there does not exist a formula in the literature to
estimate the sample size given the interval size and confidence
level which will be provided in this section. Different from the
result in \cite{per}, we begin our investigation on relative error.
Let ${\mathcal S}$ be an event, and denote $\Pr\{{\mathcal S}\}$ as
the probability associated with event ${\mathcal S}$. Specifically
we aim at solving the following problem: Let $\vep$ and $\de$ be
positive numbers less than $1$; Given margin of relative error
$\vep$ and confidence level $1 - \de$, how large should the sample
size $N$ be to ensure \be\label{pr1} \Pr \li \{  \li | \f{\wh{\mu}_N
- \mu}{\mu} \ri | < \vep   \ri \}
>1 - \de?
\ee
Our first result is the explicit formula for computing the sample
size as presented next.

\vspace{1mm}

\beT \la{main}
Let $\vep$ and $\de$ be positive numbers less than
$1$. Then the probability inequality $(\ref{pr1})$ holds, if
\be \la{bound}
N > \f{ 2(3 + 2 \vep) \ln \f{2}{\de}}
{ 3 \vep^2  }
\ee
\eeT

\vspace{1mm}

\noindent {\bf Remark 1} \
Theorem \ref{main} shows that for small values of $\vep$ and
$\de$, the {\it a priori} bound on the sample size
in (\ref{bound}) can be approximated as $\f{ 2
\ln \f{2}{\de}   } { \vep^2 }$. As a result, the sample size is
roughly inversely proportional to the square of the margin of
relative error and is linear with respect to the logarithm of the
inverse of the gap between the confidence level and $1$.  In
addition, an important feature of Theorem \ref{main} is that the sample size
bound is
independent of the value of the true second moment.

\vspace{3mm}

Theorem \ref{main} does not use the margin of absolute error as the
measure of accuracy. Its reason lies in the fact that if the margin
of absolute error is used then the sample size is essentially
dependent on the unknown second moment $\mu$, which is to be
estimated. Quite contrary, when the margin of relative error is
used, the sample size can be determined without any knowledge of the
unknown $\mu$. Moreover the relative error is a better indicator of
the quality of the estimator: an estimate with a small absolute
error may not be acceptable if the true value is also small.
Nevertheless Theorem \ref{main} will be used to provide an interval
estimate for $\mu$ given the sample size and the confidence level or
provide sample complexity bound given the estimation interval and
confidence level that is similar to a result in \cite{per} but
without complicated computation. The sample complexity result for
the interval estimate will be presented at the end of this section
after the proof for Theorem \ref{main}. First we will state and
prove the following result that indicates that the {\it a priori}
bound in (\ref{bound}) is tight as $\vep$ and $\de$ become small.

\begin{theorem} \la{main2} Under the
same hypothesis of {\rm Theorem} $\ref{main}$, define $\udl{N}(\vep,
\de)$ as an integer such that \be \udl{N}(\vep, \de) := \min  \li \{
N>0 : \; \; \; \Pr  \li \{  \li | \f{\wh{\mu}_N - \mu}{\mu}  \ri | <
\vep   \ri \} \geq 1 - \de \ri \} \ee Let $\ovl{N} (\vep, \de) = \f{
2(3 + 2 \vep) \ln \f{2}{\de} } { 3 \vep^2  }$. Then for each $\mu
> 0$ and with $Z_\de$ satisfying $\frac{1}{\sqrt{2\pi}}
\int_{-\infty}^{Z_\de} e^{-x^2/2} \ dx = 1 - \f{\de}{2}$,
\be\label{tmp1} \lim_{\de \to 0}\lim_{\vep \to 0 } \f{ \ovl{N}
(\vep, \de) } { \udl{N} (\vep, \de) } = \lim_{\de \to 0} \f{ 2 \ln
\f{2}{\de} } { Z_\de^2 } = 1 \ee
\end{theorem}

\vspace{1mm}

\noindent Proof. By the central limit theorem, we can write
\[
\Pr \li \{ \li |
\f{\wh{\mu}_N - \mu}{\mu} \ri | < \vep \ri \} =  1 -
\f{2}{\sq{2\pi}} \int_{\vep\sqrt{N}}^\iy \exp \li ( -\f{x^2}{2}
\ri ) dx + \ze(N)
\]
where $\ze(N)$ is a function of $N$ such that
$\lim_{N \to \iy} \ze(N) = 0$. Recall $Z_\de$ in the
statement of the theorem.
It follows that the minimum sample
size $\udl{N} = \udl{N} (\vep, \de)$ satisfies
\[
\f{1}{\sq{2\pi}}
\int_{\vep\sqrt{\udl{N}}}^\iy \exp \li ( -\f{x^2}{2} \ri ) dx =
\f{ \de + \ze(\udl{N}) } { 2 }
\; \; \; \; \Longrightarrow \; \; \; \;
\udl{N} (\vep,
\de) = \f{  Z_{\de + \ze( \udl{N} ) }^2 } {\vep^2}
\]
Since $Z_\de$ is a continuous function of
$\de$ and $\ze(\udl{N})
\to 0$ as $\vep \to 0$ or $N\rightarrow\infty$, we
have
\[
\lim_{\vep \to 0} \f{ \ovl{N} (\vep, \de) } { \udl{N} (\vep, \de)
} = \lim_{\vep \to 0} \f{ 2 ( 1 + \f{2}{3} \vep) \ln \f{2}{\de} }
{ Z_{\de + \ze( \udl{N} ) }^2  } = \f{ 2 \ln \f{2}{\de} } {
Z_\de^2 }
\]
that verifies the first half of (\ref{tmp1}). Alternatively since
\[
\f{1}{\sq{2 \pi}} \int_{Z_\de}^\iy \exp \li ( - \f{x^2}{2} \ri  )
d x = \f{\de}{2}
\]
we have that $Z_\de \to \iy$ as $\de \to 0$.
Consequently,
\[
\frac{\delta}{2} = \f{1}{\sq{2 \pi}} \int_{Z_\de}^\iy \exp
\li ( - \f{x^2}{2} \ri ) d x \; \; \; \rightarrow \; \; \;
\f{1}{\sq{2 \pi} Z_\de} \int_{Z_\de}^\iy \exp \li(-\frac{x^2}{2}\ri)
d\li(\frac{x^2}{2}\ri) =
\f{1}{\sq{2 \pi} Z_\de} \exp \li
( - \f{Z_\de^2}{2} \ri )
\]
as $\de \to 0$. Taking the natural logarithm on both sides with
appropriate rearrangement yields
\[
\lim_{\de \to 0} \f{ 2\ln \f{2}{\de} } {Z_\de^2} = \lim_{\de \to 0}
\li [ 1 + \f{\ln (2\pi Z_\de^2 )}{ Z_\de^2 }  \ri ] =
1
\]
and therefore (\ref{tmp1}) follows that concludes the
proof. \pfbox

\vspace{3mm}

\noindent {Proof of Theorem \ref{main}:} \ It can be seen that

\be\label{upb} \Pr  \left\{\left| \f{\wh{\mu}_N - \mu}{\mu}\right|
\geq \vep  \right\}  = \Pr \left\{ \sum_{i = 1}^N \f{X_i^2}{\mu}
\geq N (1 + \vep)\right\} + \Pr \left\{ \sum_{i = 1}^N
\f{X_i^2}{\mu} \leq N (1 - \vep)\right\} \ee
 \noindent Note that
$\f{2X_i^2}{\mu}$ has mean 1 and possesses a $\chi^2$ distribution
with $2$ degrees of freedom. Because $\left\{ \f{2 X_i^2}{\mu}
\right\}_{i = 1}^N$ are i.i.d., random variable $R = \sum_{i = 1}^N
\f{2 X_i^2}{\mu}$ possesses a $\chi^2$ distribution with $2N$
degrees of freedom. By equation (2.1-138) in page 46 of
\cite{proakis}, we have
\[
\Pr \li
\{ \sum_{i = 1}^N X_i^2 \geq N (1 + \vep) \mu \ri \}  = e^{ - N(1
+ \vep)} \sum_{k = 0}^{N-1} \f{ [N(1 + \vep)]^k } { k!}
\]
Define a Poisson random variable $Y$ with mean $\se = N(1 + \vep)$.
Then
\[
\Pr \li
\{ \sum_{i = 1}^N X_i^2 \geq N (1 + \vep) \mu \ri \} = \Pr \{ Y < N \}
\leq  \inf_{\lm < 0} \bb{E} [ e^{ \lm (Y - N) } ]  =  \inf_{\lm <
0} e^{- \se} e^{ \se e^\lm - \lm N }
\]
where the Chernoff bound from \cite{chernoff} is applied.
It is easy to see that the infimum is
achieved at $\lm = \ln ( \f{N}{\se} ) < 0$ for which
we have $e^{- \se} e^{ \se e^\lm - \lm N } = e^{- \se}\li(
\f{ \se e } { N }\ri)^N$. Consequently $\Pr \{ Y < N \} \leq e^{-
\se} \li( \f{ \se e } { N } \ri)^N$ and
\be
\Pr \li\{ \sum_{i = 1}^N
X_i^2 \geq N (1 + \vep)\mu  \ri \}
=\Pr \{ Y < N \}
\leq e^{- \se} \li( \f{\se e}{N}\ri)^N
= \li[ (1+ \vep) e^{- \vep} \ri]^N < \f{\de}{2}
\ee
provided that $N > \f{ \ln (\f{2}{\de})   } { \vep - \ln(1 + \vep)  }$.

\vspace{1mm}

Similarly define a Poisson random variable $Z$ with mean $\vse
= N(1 - \vep)$. Then the Chernoff bound leads to
$\Pr \{ Z \geq N \} \leq e^{- \vse} \li ( \f{
\vse e } { N } \ri )^N$ as earlier. It follows that
\be
\Pr \li \{
\sum_{i = 1}^N X_i^2 \leq N (1 - \vep) \mu  \ri \}=
\Pr \{ Z \geq N \} \leq e^{- \vse} \li ( \f{\vse e}{N}
\ri )^N = e^{- \vse} \li [ (1 - \vep) e \ri ]^N=
\li [ (1- \vep) e^{\vep} \ri ]^N < \f{\de}{2}
\ee
provided that $N
 > \f{ \ln (\f{2}{\de})   } { - \vep - \ln(1 - \vep)  }$.
Using the following inequalities (to be proven later)
\be
\la{inq} \f{ 3\vep^2 } { 2(3 + 2\vep) } < \vep - \ln(1 +
\vep) < - \vep - \ln(1 - \vep) \; \; \; \forall \; \vep\ \in (0, \; 1)
\ee
we have that if $N > \f{ 2(3 +
2 \vep) \ln (\f{2}{\de})  } { 3 \vep^2  }$, then
$N
 > \f{ \ln (\f{2}{\de})   } { \vep - \ln(1 + \vep)  }
>\f{ \ln (\f{2}{\de})   } { - \vep - \ln(1 - \vep)  }$, and consequently
\[
\Pr \li\{ \sum_{i =
1}^N X_i^2 \geq N (1 + \vep) \mu \ri\} < \f{\de}{2},
\; \; \; \; \; \Pr \li\{ \sum_{i = 1}^N X_i^2 \leq N (1 - \vep) \mu  \ri\} <
\f{\de}{2}
\]
leading to $\Pr \li\{\li| \f{\wh{\mu}_N - \mu}{\mu} \ri| <
\vep \ri\} > 1 - \de$.  Finally to prove the first inequality of
(\ref{inq}), define
\[
g(\vep) = \vep - \ln(1 + \vep) - \f{ 3\vep^2
} { 2(3+2\vep) }
\]
Then $g(0) = 0$ and
\[
\f{d g(\vep)} { d \vep} = \f{ 10\vep^2  } { (1 + \vep)  (3
+ 2\vep)^2}
>0
\]
It follows that the first inequality of (\ref{inq}) holds.
Similarly to prove the second inequality of (\ref{inq}), define
\[
h(\vep) = \vep - \ln(1 + \vep) - [  - \vep - \ln(1 - \vep)]
\]
Then $h(0) = 0$ and
\[
\f{d h(\vep)} { d \vep} = - \f{ 2 \vep^2  } { 1 - \vep^2 } < 0
\; \; \; \forall \; \vep\ \in (0, \; 1)
\]
Hence the second inequality of (\ref{inq}) is true.
The proof of Theorem \ref{main} is now complete. \pfbox

\vspace{3mm}

\noindent {\bf Remark 2} \ The {\it a priori} sample complexity bound in
Theorem \ref{main}, while being tight asymptotically, can be further
improved by taking $N=N_{\rm m1}$ to be the minimum positive integer
such that
\be\label{nmin1} F_{\chi^2, 2N_{\rm m1}} ( (1+ \vep)
2N_{\rm m1}) - F_{\chi^2, 2N_{\rm m1}} ( (1 - \vep) 2N_{\rm m1}) > 1
- \de
\ee
in light of (\ref{upb}) where $F_{\chi^2, 2N}(\cdot)$ is
the CDF of $\chi^2$ random variable with $2N$ degree of freedom.
Indeed there holds $N> \f{ 2(3 + 2 \vep) \ln(\f{2}{\de})} { 3
\vep^2  }> \f{\ln(\f{2}{\de})}{\vep-\ln (1+\vep)} \geq N_{\rm m1}$. In addition for
positive integer $N<N_{\rm m1}$, \be \Pr \li \{  \li | \f{\wh{\mu}_N
- \mu}{\mu} \ri | < \vep   \ri \} \leq 1-\de \ee

\vspace{3mm}

Although only the margin of relative error is investigated in
Theorem \ref{main} for the sample complexity problem in statistical
estimation of the 2nd moment of Rayleigh random variables, a similar
result can be obtained as well for the corresponding interval
estimate. The probability associated with
such an interval estimate will be termed as {\it coverage
probability} to signify its utility in coverage analysis \cite{b, c}.

\begin{corollary} \label{co1}
{\rm (i)} Let $\vep$ and $\de$ be positive numbers less than
$1$. Then
\be\label{tmp2}
\Pr \li \{ \frac{\wh\mu_N}{1+\vep}< \mu <
\frac{\wh\mu_N}{1-\vep} \ri \}
>1 - \de
\ee provided that $N$ satisfies the sample complexity bound in
$(\ref{bound})$. Recall that $\widehat{\mu}_N$ as in $(\ref{muhat})$
is the maximum likelihood estimate for the $2$nd moment.

\vspace{2mm}

\noindent {\rm (ii)} \ Let $\de \in (0,1)$, $N > \f{10}{3} \ln
(\f{2}{\de})$, and
\be\label{ert}
\vep = \f{ 2 \ln \f{2}{\de} } { 3 N } \li ( 1 + \sq{1 + \f{ 9 N }
{2 \ln \f{2}{\de} } } \ri )
\ee
Then the inequality $(\ref{tmp2})$ for the coverage probability holds.
\end{corollary}

\noindent Proof. We note the following chain of equalities:
\bee
\Pr \li \{\li
| \f{\wh{\mu}_N - \mu}{\mu} \ri | < \vep   \ri \}
&=& \Pr \li \{ -\vep< \f{\wh{\mu}_N - \mu}{\mu}< \vep \ri\}  \\
&=& \Pr \li \{ (1-\vep)\mu< \wh{\mu}_N < (1+\vep)\mu \ri \}
=  \Pr \li \{ \frac{\wh\mu_N}{1+\vep} < \mu <
\frac{\wh\mu_N}{1-\vep} \ri \}
\eee
Hence (i) holds true in light of Theorem \ref{main}.
For (ii) solving equation $N = \f{ 2(3
+ 2 \vep) \ln(\f{2}{\de})} { 3 \vep^2  }$ with respect to $\vep$
yields the unique positive root $\vep$ as in (\ref{ert}).
Note that $\f{ 2(3 + 2 \vep) \ln(\f{2}{\de})} { 3 \vep^2  }$
decreases as $\vep > 0$ increases, and $\vep=1$ is the unique positive root
of $N = \f{2(3 + 2 \vep)\ln(\f{2}{\de})}{3\vep^2}$
at $N=\f{10}{3} \ln(\f{2}{\de})$. Therefore, the root $\vep$ is
less than $1$ for $N
>\f{10}{3} \ln(\f{2}{\de})$. It follows that (ii) holds true as well
in light of Theorem \ref{main}. \pfbox

\vspace{3mm}

\noindent {\bf Remark 3} \

The problem of interval estimate has been investigated in \cite{per}
that proposes a computational procedure but without an {\it a
priori} bound on the sample size as in Corollary \ref{co1}. Indeed
it computes the exact value of $\Pr \li \{ C_1\wh\mu_N\leq \mu\leq
C_2\wh\mu_N \ri \} = p$ given $0<p<1$
by searching for $(C_1,
C_2)$ which
are dependent on $N$ and $\vep$ and satisfy
$0<C_1<1<C_2$. Although in a different form, $C_1$ and
$C_2$ approach to 1 in accordance with Theorem \ref{main2}
as $N\rightarrow\infty$. However
the algorithm as proposed in \cite{per} makes the use of the tabulated
$\chi^2$ CDF repeatedly in searching for an appropriate interval
and sample size that admit multiple solutions for a given $p$
because of the trade-off between the interval length and sample size.
Such a search is
bypassed in Corollary \ref{co1}.
In fact such a search is
unnecessary that is completely characterized by (\ref{nmin1}).
The explicit relation between $N$
and $(\vep, \de)$ enables us to compute an interval estimate for
$\mu$ {\it a priori} given the sample size and the confidence level,
or compute the {\it a priori} sample complexity bound given the
estimation interval and confidence level that is in sharp contrast
to the existing results in the literature. Numerical result
indicates that our explicit formula in Corollary \ref{co1} on the
sample complexity bound is very tight as shown on the left of Figure 1
with $\de=0.01$ where $C_1$ and $C_2$ are solved from (\ref{nmin1}) and
\[
\f{1}{(1+\vep)C_1} = \f{\f{\wh{\mu}_N}{1+\vep}}{C_1\wh{\mu}_N},
\; \; \; \; \;
\f{1}{(1-\vep)C_2} = \f{\f{\wh{\mu}_N}{1-\vep}}{C_2\wh{\mu}_N}
\]
On the right of Figure 1 plots $\f{1+\vep}{1-\vep}$ in unit dB
vs. the {\it a priori} bound $N$ that is a function of
$\vep$, given $\delta$.

\vspace*{2.1in}
\begin{figure}[ht]\label{fig0}
\includegraphics{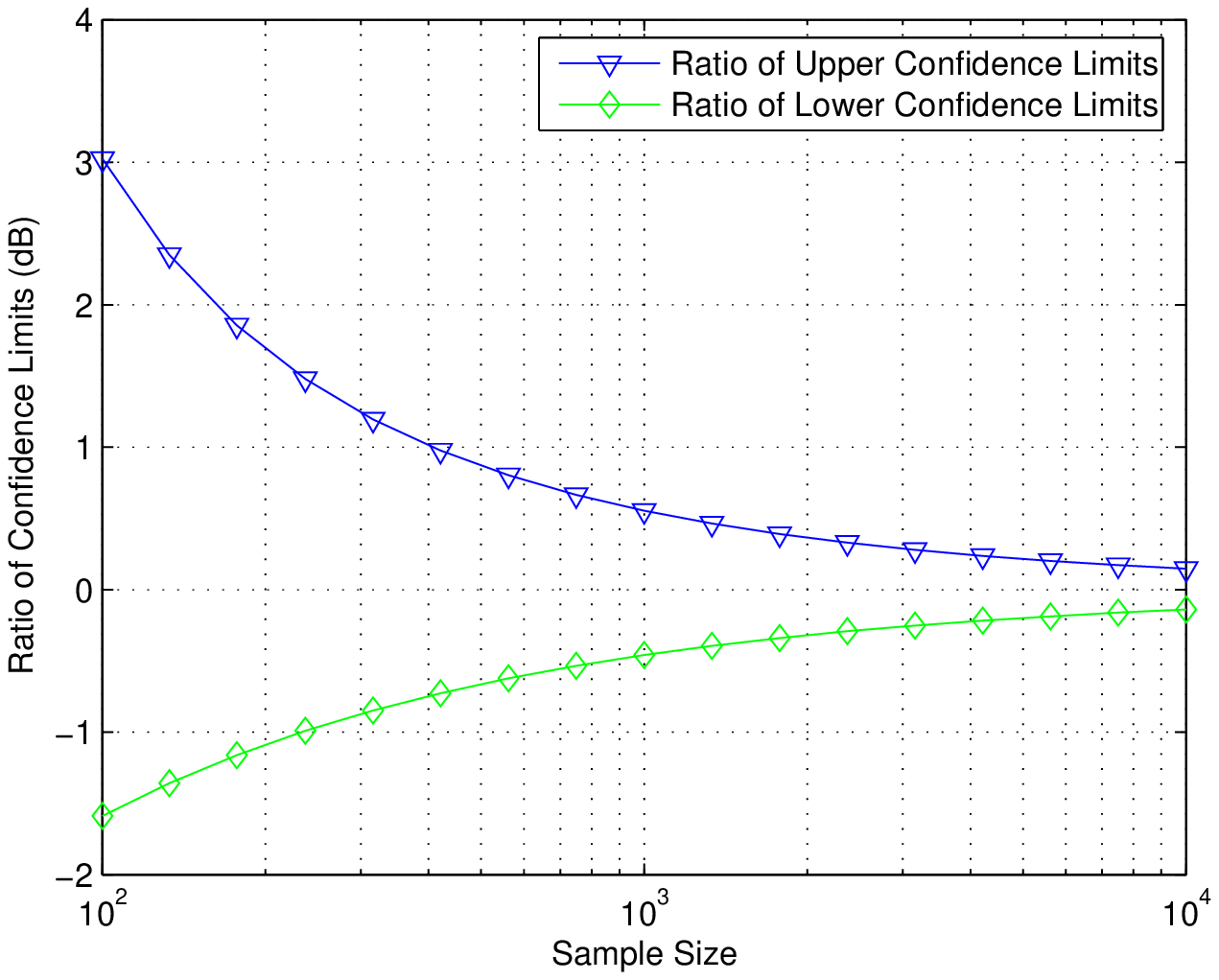}
\includegraphics{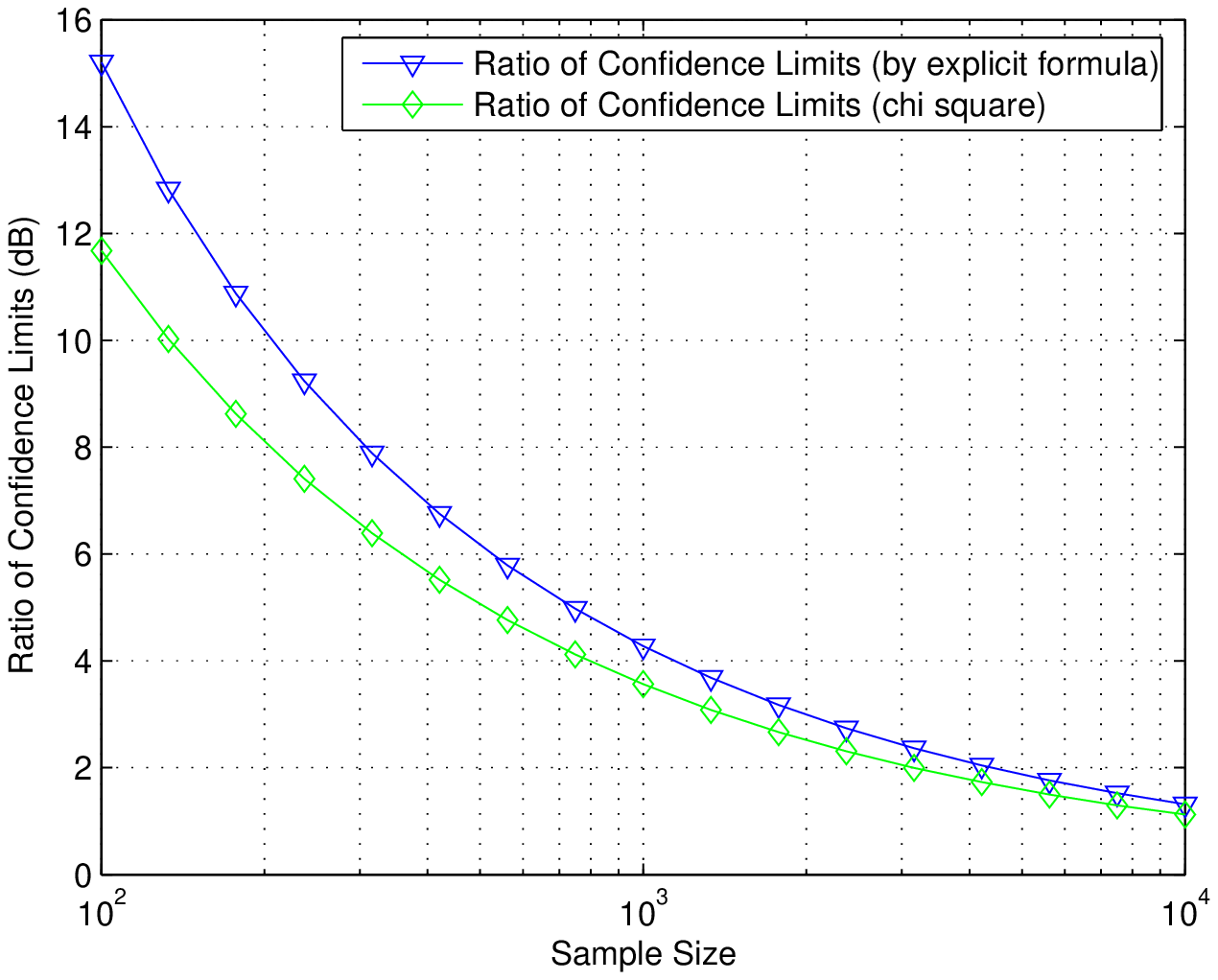}
\caption{Plots of $20\log_{10}\left(\f{1}{(1+\vep)C_1}\ri)$ and
$20\log_{10}\left(\f{1}{(1-\vep)C_2}\ri)$ (left), and
$20\log_{10}\left(\f{1+\vep}{1-\vep}\right)$ (right)
vs. $N$}
\end{figure}

\noindent In summary our results on the sample complexity
are quite satisfactory as demonstrated in Figure 1
that can in fact be extended to the noisy case.

\section{Measurement and Estimation with Noisy Samples}

Recall that in the noisy case two sets of i.i.d.
BPSK symbols are transmitted with difference in their average power
satisfying $P_{s1}>P_{s2}\geq 0$. The problem is statistical estimation for
$\si_H^2$, $\si_V^2$, and the SNR, as well as their associated sample
complexity bound based on i.i.d. observations of the received signal
$X = HS + V$ which are recorded as $\{X_{1,i}\}_{i=1}^{N_1}$
and $\{X_{2,i}\}_{i=1}^{N_2}$, respectively. We consider
the case where $N_1=N_2=N$ with $\overline{X_{1,N}^2}$
and $\overline{X_{2,N}^2}$ as defined in (\ref{tmp0}). Two
different subjects will be investigated.

\subsection{Interval Estimate}

We will demonstrate not only how to construct a simultaneous
confidence interval based on the noisy measurements but also how
tight these interval estimates are. Our main results begin with the
following theorem that provides a simultaneous interval estimate and
the solution to the sample complexity problem with an {\it a priori}
bound.

\beT \la{Th2} Let $\vep$ and $\de$ be positive numbers less than
$1$. Denote $\Delta P_s=P_{s1}-P_{s2}$. Then
\be\label{pr2}
\Pr \li \{ \begin{array}{c} \hspace*{6mm}
\f{1}{\Delta P_s} \li(\f{\overline{X_{1,N}^2}}
{1 + \vep} - \f{\overline{X_{2,N}^2}}{1 - \vep} \ri)
< \si_H^2 < \f{1}{\Delta P_s} \li
(\f{\overline{X_{1,N}^2}}{1 - \vep} -
\f{\overline{X_{2,N}^2}}{1 + \vep} \ri) \; \, \, \& \\
\f{1}{\Delta P_s} \li
(\f{P_{s1}\overline{X_{2,N}^2}}{1 + \vep} -
\f{P_{s2}\overline{X_{1,N}^2}}{1 - \vep} \ri)
< \si_V^2 < \f{1}{\Delta P_s} \li
(\f{P_{s1}\overline{X_{2,N}^2}}{1 - \vep} -
\f{P_{s2}\overline{X_{1,N}^2}}{1 + \vep} \ri)
\end{array}\ri \}
>1 - \de
\ee
provided that $N
>\f{2(3 + 2 \vep) \ln (\f{4}{\de})} {3\vep^2}$.
\eeT

\noindent Proof. For BPSK the transmitted symbol $S$ is real
with the transmitting power $P_s=|S|^2$, and thus
\[
X = HS + V = (H_RS + V_R) + j (H_IS + V_I) = \sq{(P_s\si_H^2 +
\si_V^2)} (A+jB)
\]
where $A$ and $B$ are independent
Gaussian variables with zero
mean and unity variance. Hence
$|X|^2 = (P_s\si_H^2 +
\si_V^2)(A^2 + B^2)$ is Rayleigh distributed.
It follows that $\{X_{1, i}\}_{i=1}^N$ and $\{X_{2, i}\}_{i=1}^N$
are both i.i.d. observations of $X$, and with $\overline{X_{1,N}^2}$
and $\overline{X_{2,N}^2}$ as defined in (\ref{tmp0}),
\[
\bb{E} \li [\overline{X_{\ell,N}^2}\ri] =
\left(P_{s_\ell}\si_H^2 + \si_V^2\right),
\; \; \; \; \; \ell=1, 2
\]
Applying (i) of Corollary \ref{co1} yields the interval estimate
\[
\Pr \li \{ \f{\overline{X_{\ell,N}^2}}{1 + \vep} <
(P_{s_\ell}\si_H^2 + \si_V^2) < \f{\overline{X_{\ell,N}^2}}{1 - \vep}
\ri \} > \sq{1 - \de}, \; \; \; \; \; \ell=1, 2
\]
for an appropriately chosen sample size $N$.  In fact the
above holds, if
\be\label{tmp4}
N > \frac{2(3+2\vep){\rm
ln}\left(\frac{2}{1-\sqrt{1-\delta}}\right)}{3\vep^2}
\ee
in light of Corollary \ref{co1}. The fact that
$(1-\sqrt{1-\de})^2=2-\de-2\sqrt{1-\de}>0$ implies that
$1-\sqrt{1-\delta} > \frac{\delta}{2}$. Hence the inequality
in (\ref{tmp4}) holds provided that $N>
\f{2(3 + 2 \vep) \ln (\f{4}{\de})} {3\vep^2}$ as in the
statement of the theorem. Now by the independence of the experiment,
\[
\Pr \li  \{
\f{\overline{X_{1,N}^2}}{1 + \vep} <
(P_{s1}\si_H^2 + \si_V^2) < \f{\overline{X_{1,N}^2}}{1 - \vep},
\; \; \; \f{\overline{X_{2,N}^2}}{1 + \vep} <
(P_{s2}\si_H^2 + \si_V^2) < \f{\overline{X_{2,N}^2}}{1 - \vep}\ri\}
>1 - \de
\]
Define the following sets of events as
illustrated in Figure 2:
\begin{eqnarray}\label{s0}
{\mathcal S}_0 &=& \left\{(\sigma_H^2, \sigma_V^2): \; \;
\mbox{the joint events in the probability expression (\ref{pr2})
hold}\right\}\\
\label{s1} {\mathcal S}_1 &=& \left\{(\sigma_H^2, \sigma_V^2): \; \;
\f{\overline{X_{1,N}^2}}{1 + \vep} <
(P_{s1}\si_H^2 + \si_V^2) < \f{\overline{X_{1,N}^2}}{1 - \vep}\right\}\\
\label{s2} {\mathcal S}_2 &=& \left\{(\sigma_H^2, \sigma_V^2): \; \;
\f{\overline{X_{2,N}^2}}{1 + \vep} <
(P_{s2}\si_H^2 + \si_V^2) <
\f{\overline{X_{2,N}^2}}{1 - \vep}\right\}
\end{eqnarray}
We will show that ${\mathcal S}_1\cap {\mathcal S}_2\subseteq {\mathcal
S}_0$
leading to $\Pr\{{\mathcal S}_0\}\geq \Pr\{{\mathcal S}_1\cap {\mathcal
S}_2\}
>1-\de$ thereby concluding the proof.

\vspace*{1.6in}
\begin{figure}[ht]\label{fig1}
\includegraphics{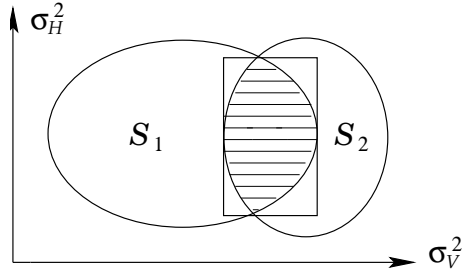}
\caption{Illustration of the sets of events}
\end{figure}

\vspace{2mm}

\noindent
Indeed suppose that $(\sigma_H^2, \sigma_V^2) \in {\mathcal S}_1\cap
{\mathcal S}_2$
represented by the shaded area in Figure 2. Then
\[
\f{\overline{X_{1,N}^2}}{1 + \vep} <
(P_{s1}\si_H^2 + \si_V^2) < \f{\overline{X_{1,N}^2}}{1 - \vep}
\; \; \; \& \; \; \; -\f{\overline{X_{2,N}^2}}{1 - \vep} <
-(P_{s2}\si_H^2 + \si_V^2) < -\f{\overline{X_{2,N}^2}}{1 + \vep}
\]
The above inequalities can be combined to eliminate
$\sigma_V^2$ yielding
\be\label{sh}
(\sigma_H^2, \sigma_V^2) \in {\mathcal S}_H:=\left\{\sigma_H^2: \; \; \;
\f{1}{\Delta P_s}\left(\f{\overline{X_{1,N}^2}}
{1 + \vep} - \f{\overline{X_{2,N}^2}}{1 - \vep}\ri)
< \si_H^2 < \f{1}{\Delta P_s}
\left(\f{\overline{X_{1,N}^2}}{1 - \vep} -
\f{\overline{X_{2,N}^2}}{1 + \vep}\ri)
\right\}
\ee
On the other hand for $(\sigma_H^2, \sigma_V^2)
\in {\mathcal S}_1\cap {\mathcal S}_2$, there hold
\bee
\f{P_{s1}\overline{X_{2,N}^2}}{1 + \vep} <
&(P_{s1}P_{s2}\si_H^2 + P_{s1}\si_V^2)& <
\f{P_{s1}\overline{X_{2,N}^2}}{1 - \vep}
\; \; \; \& \\
-\f{P_{s2}\overline{X_{1,N}^2}}{1 - \vep} <
&
-(P_{s1}P_{s2}\si_H^2 + P_{s2}\si_V^2)& <
-\f{P_{s2}\overline{X_{1,N}^2}}{1 + \vep}
\eee
The above inequalities can be combined to eliminate
$\sigma_H^2$ yielding
$(\sigma_H^2, \sigma_V^2) \in {\mathcal S}_V$ with
\be\label{sv}
{\mathcal S}_V:= \left\{\sigma_V^2: \; \; \;
\f{1}{\Delta P_s} \li
(\f{P_{s1}\overline{X_{2,N}^2}}{1 + \vep} -
\f{P_{s2}\overline{X_{1,N}^2}}{1 - \vep} \ri)
< \si_V^2 < \f{1}{\Delta P_s} \li
(\f{P_{s1}\overline{X_{2,N}^2}}{1 - \vep} -
\f{P_{s2}\overline{X_{1,N}^2}}{1 + \vep} \ri)\right\}
\ee
We can thus conclude that ${\mathcal S}_1\cap {\mathcal S}_2\subseteq
{\mathcal S}_0 =
{\mathcal S}_H\cap{\mathcal S}_V$ represented by the rectangular area
in Figure 2. The proof is now complete.
\pfbox

\vspace{3mm}

The proof of Theorem \ref{Th2} shows that two signal sets of
different strengths are necessary in order to resolve two different
interval estimates for $\sigma_H^2$ and $\sigma_V^2$, respectively.
Clearly signals of more than two different strengths
become necessary in order to acquire more
than two independent parameters in
statistical estimation of wireless channels.
Theorem \ref{Th2} also shows that
for
$N>\f{2(3+\vep)\ln(\f{4}{\de})}{3\vep^2}$,
\bee
\Pr \li \{{\mathcal S}_H\ri\} &=&
\Pr \li \{\f{1}{\Delta P_s} \li(\f{\overline{X_{1,N}^2}}
{1 + \vep} - \f{\overline{X_{2,N}^2}}{1 - \vep} \ri)
< \si_H^2 < \f{1}{\Delta P_s} \li
(\f{\overline{X_{1,N}^2}}{1 - \vep} -
\f{\overline{X_{2,N}^2}}{1 + \vep} \ri)\ri\} > 1 - \de \\
\Pr \li \{{\mathcal S}_V\ri\} &=& \Pr \li \{ \f{1}{\Delta P_s} \li
(\f{P_{s1}\overline{X_{2,N}^2}}{1 + \vep} -
\f{P_{s2}\overline{X_{1,N}^2}}{1 - \vep} \ri)
< \si_V^2 < \f{1}{\Delta P_s} \li
(\f{P_{s1}\overline{X_{2,N}^2}}{1 - \vep} -
\f{P_{s2}\overline{X_{1,N}^2}}{1 + \vep} \ri)\ri \}
>1 - \de
\eee

\vspace{1mm}

\noindent {\bf Remark 4} \ In light of Remark 2 and
the definition of $F_{\chi^2, 2N}(\cdot)$,
the sample complexity bound in Theorem \ref{Th2} can be
improved by taking $N=N_{\rm m2}$, the minimum
positive
integer such that
\be\label{nmin2}
F_{\chi^2, 2N_{\rm m2}} ( (1+ \vep)2N_{\rm m2}) -
F_{\chi^2, 2N_{\rm m2}} ( (1 - \vep) 2N_{\rm m2}) \geq \sq{1 - \de}
\ee

\vspace{1mm}

Recall that SNR $=  10 \log_{10} \bb{E}[|S|^2] + 10 \log_{10} \li (
\f{ \si_H^2 } { \si_V^2 } \ri )$, it is both important and of
independent interest to obtain an interval estimate for $\f{ \si_H^2
} { \si_V^2 }$ that is presented in the next result.

\begin{corollary} \label{co2}
Let $\vep$, $\de$ be positive numbers less than
$1$, and $P_{s2} = 0$.  Suppose that $N$ is the
minimum positive  integer such that
$(\ref{nmin2})$ holds. Then,
\[
\Pr \li \{ \f{1}{P_{s1}} \li(\f{(1 -
\vep)\overline{X_{1,N}^2}}{(1 + \vep)\overline{X_{2,N}^2}}
- 1 \ri ) < \f{ \si_H^2 } { \si_V^2 } <
\f{1}{P_{s1}} \li(\f{(1 +
\vep)\overline{X_{1,N}^2}}{(1 - \vep)\overline{X_{2,N}^2}}
- 1 \ri )\ri \} > 1-\de
\]
\end{corollary}

\noindent Proof. In reference to the sets of the events as
illustrated in Figure 2, the hypothesis $P_{s2}=0$ implies
\[
{\mathcal S}_{1}\cap{\mathcal S}_2 = \li\{(\sigma_H^2, \sigma_V^2): \; \; \;
\f{\overline{X_{1,N}^2}}{1 + \vep} <
(P_{s1}\si_H^2 + \si_V^2) < \f{\overline{X_{1,N}^2}}{1 - \vep}
\; \; \& \; \; \frac{\overline{X_{2,N}^2}}{1 + \vep} <
\si_V^2 < \f{\overline{X_{2,N}^2}}{1 - \vep}\ri\}
\]
Suppose that $(\sigma_H^2, \sigma_V^2) \in {\mathcal S}_1\cap {\mathcal
S}_2$.
Then some algebraic manipulations yield that
\be\label{shv}
(\sigma_H^2, \sigma_V^2) \in {\mathcal S}_{H,V}
:= \li\{(\sigma_H^2, \sigma_V^2): \; \; \;
\f{(1 -
\vep)\overline{X_{1,N}^2}}
{(1 + \vep)\overline{X_{2,N}^2}} <
\f{P_{s1}\si_H^2 } { \si_V^2 } + 1<
\f{(1 + \vep)\overline{X_{1,N}^2}}{(1 - \vep)\overline{X_{2,N}^2}}\ri\}
\ee
As a result ${\mathcal S}_1\cap {\mathcal S}_2\subseteq {\mathcal S}_{H,V}$
and thus $\Pr\{{\mathcal S}_{H,V}\}\geq \Pr\{{\mathcal S}_1\cap
{\mathcal S}_2\}\geq 1-\de$, if $N$ satisfies the
(\ref{nmin2}) by Remark 3. Therefore the corollary is true. \pfbox

\vspace{3mm}

In the case $P_{s2}\neq 0$, the result in Corollary
\ref{co2} needs to be replaced by a more complex expression:
\be\label{co22}
\Pr\li\{\frac{\li(1-\f{P_{s2}}{P_{s1}}\ri)\overline{X_{1,N}^2}}
{\li(\f{1+\vep}{1-\vep}\ri)P_{s1}\overline{X_{2,N}^2} -
P_{s2}\overline{X_{1,N}^2}}
- \frac{1}{P_{s1}} <\frac{\sigma_H^2}{\sigma_V^2}<
\frac{\li(1-\f{P_{s2}}{P_{s1}}\ri)\overline{X_{1,N}^2}}
{\li(\f{1-\vep}{1+\vep}\ri)P_{s1}\overline{X_{2,N}^2} -
P_{s2}\overline{X_{1,N}^2}}
- \frac{1}{P_{s1}}
\ri\} > 1-\de
\ee
if the positive integer $N$ satisfies (\ref{nmin2}),
and if $\li(\f{1+\vep}{1-\vep}\ri)P_{s1}\overline{X_{2,N}^2}
- P_{s2}\overline{X_{1,N}^2}>0$ and
$\li(\f{1-\vep}{1+\vep}\ri)P_{s1}\overline{X_{2,N}^2}
- P_{s2}\overline{X_{1,N}^2}>0$.
The above reduces to that of
Corollary \ref{co2} in the case $P_{s2}=0$.

\vspace{1mm}

We investigate next whether or not the interval
estimates presented in Theorem \ref{Th2} and Corollary \ref{co2}
are tight. Such an issue is clearly dependent on the average
transmitting powers $P_{s1}$ and $P_{s2}$. Recall that
$P_{s1}>P_{s2}\geq 0$. Intuitively the larger
$P_{s1}$ and the smaller $P_{s2}$ are, the
tighter the interval estimates are that is validated
by the following result.

\beT \la{Th3} Let $\vep$, $\de$ be positive numbers less than
$1$, and $N = N_{\rm m2}-1\geq N_{\rm m1}$ with
$N_{\rm m1}$ and $N_{\rm m2}$ the minimum positive integers
satisfying $(\ref{nmin1})$ and $(\ref{nmin2})$, respectively. Then
\be\la{ineq1}
1-\de \leq \Pr \li \{\f{1}{\Delta P_s} \li(\f{\overline{X_{1,N}^2}}
{1 + \vep} - \f{\overline{X_{2,N}^2}}{1 - \vep} \ri)
< \si_H^2 < \f{1}{\Delta P_s} \li
(\f{\overline{X_{1,N}^2}}{1 - \vep} -
\f{\overline{X_{2,N}^2}}{1 + \vep} \ri)\ri\} <1 - \frac{\de}{2}
\ee
for sufficiently large $P_{s1}$ with constant $P_{s2}$. In addition
\be\la{ineq2}
1-\de \leq \Pr \li \{ \f{1}{\Delta P_s} \li
(\f{P_{s1}\overline{X_{2,N}^2}}{1 + \vep} -
\f{P_{s2}\overline{X_{1,N}^2}}{1 - \vep} \ri)
< \si_V^2 < \f{1}{\Delta P_s} \li
(\f{P_{s1}\overline{X_{2,N}^2}}{1 - \vep} -
\f{P_{s2}\overline{X_{1,N}^2}}{1 + \vep} \ri)\ri \}
< 1 - \frac{\de}{2}
\ee
for sufficiently small $P_{s2}$ with constant $P_{s1}$.
Finally $1-\f{\de}{2}$ is also an upper bound for
the confidence level in {\rm Theorem} $\ref{Th2}$
and {\rm Corollary} $\ref{co2}$ with sufficiently large
$P_{s1}$ and sufficiently small $P_{s2}$, if $N=N_{\rm m2}-1$. \eeT

\noindent Proof. In light of Remark 4 and the hypothesis on $N$,
there holds \bee \sqrt{1-\de} &\geq &
\Pr\li\{\f{\overline{X_{1,N}^2}}{1+\vep}
< (P_{s1}\sigma_H^2+\sigma_V^2) < \f{\overline{X_{1,N}^2}}{1-\vep}\ri\}\\
&=& \Pr\li\{(1-\vep)(P_{s1}\sigma_H^2+\sigma_V^2) <
\overline{X_{1,N}^2}< (1+\vep)(P_{s1}\sigma_H^2+\sigma_V^2)\ri\} \\
&=& \Pr\li\{(1-\vep)\li(\sigma_H^2+\frac{\sigma_V^2}{P_{s1}}\ri) <
\f{\overline{X_{1,N}^2}}{P_{s1}} <
(1+\vep)\li(\sigma_H^2+\frac{\sigma_V^2}{P_{s1}}\ri)\ri\} \\
&\rightarrow& \Pr\li\{(1-\vep)\sigma_H^2 < \wh{\sigma}_{H,N}^{2,(\infty)} <
(1+\vep)\sigma_H^2\ri\}
>1-\de
\eee where $\wh{\sigma}_{H,N}^{2,(\infty)}$ is the limit of
$\wh{\sigma}_{H,N}^2$, a point estimate for $\sigma^2_H$ to be
investigated later: \be \wh{\sigma}_{H,N}^{2,(\infty)} =
\lim_{P_{s1}\rightarrow\infty} \f{\overline{X_{1,N}^2}}{P_{s1}} =
\lim_{P_{s1}\rightarrow\infty} \li\{ \wh{\sigma}_{H,N}^2 :=
\f{\overline{X_{1,N}^2}- \overline{X_{2,N}^2}}{P_{s1}-P_{s2}}\ri\}
\ee The last inequality follows from the hypothesis on $N_{\rm
m2}>N_{\rm m1}$ and from the fact that $P_{s1}\rightarrow\infty$
corresponds to the noise-free case in Section 3. Recall ${\mathcal
S}_H$ as defined in (\ref{sh}) that can be alternatively described as
\be {\mathcal S}_H = \li\{\sigma_H^2: \; \; \de
P_s(1-\vep)\sigma_H^2 +
\f{(1-\vep)\overline{X_{2,N}^2}}{(1+\vep)P_{s1}} <
\f{\overline{X_{1,N}^2}}{P_{s1}}< \de P_s(1+\vep)\sigma_H^2 +
\f{(1+\vep)\overline{X_{2,N}^2}}{(1-\vep)P_{s1}}\ri\} \ee where $\de
P_s = \li(1-\f{P_{s2}}{P_{s1}}\ri)$. It is easy to see that as
$P_{s1}\rightarrow\infty$,
\[
\Pr\li\{{\mathcal S}_H\ri\} \rightarrow
\Pr\li\{(1-\vep)\sigma_H^2 < \wh{\sigma}_{H,N}^{2,(\infty)} <
(1+\vep)\sigma_H^2\ri\} > 1-\delta
\]
It follows that the lower bound in
(\ref{ineq1}) is true by taking $P_{s1}$
sufficiently large. On the other hand
\[
\Pr\li\{{\mathcal S}_H\ri\} \rightarrow
\Pr\li\{(1-\vep)\sigma_H^2 < \wh{\sigma}_{H,N}^{2,(\infty)} <
(1+\vep)\sigma_H^2\ri\} \leq \sqrt{1-\de}
\]
It follows that there exists $P_{s1}>0$ such that
\[
\Pr\li\{{\mathcal S}_H\ri\} < \sqrt{1-\de}
+ \frac{1}{2}\li(1-\sqrt{1-\de}\ri)^2 = 1-\frac{\de}{2}
\]
provided that $P_{s1}$ is large enough
that concludes the upper bound in (\ref{ineq1}).
The proof for (\ref{ineq2}) is similar. Indeed
by applying the result from the previous section and
using the hypothesis $N_{\rm m1}<N_{\rm m2}$,
\bee
\sqrt{1-\de} &\geq & \Pr\li\{\f{\overline{X_{2,N}^2}}{1+\vep}
< (P_{s2}\sigma_H^2+\sigma_V^2) < \f{\overline{X_{2,N}^2}}{1-\vep}\ri\}\\
&=& \Pr\li\{(1-\vep)(P_{s2}\sigma_H^2+\sigma_V^2) <
\overline{X_{2,N}^2}< (1+\vep)(P_{s2}\sigma_H^2+\sigma_V^2)\ri\} \\
&\rightarrow& \Pr\li\{(1-\vep)\sigma_V^2 <
\overline{X_{2,N}^2} <
(1+\vep)\sigma_V^2\ri\}
>1-\de
\eee
as $P_{s2}\rightarrow 0$. Recall ${\mathcal S}_V$
as defined in (\ref{sv}). Then
\[
1-\f{\de}{2} > \sqrt{1-\de} \geq \lim_{P_{s2}\rightarrow 0}
\Pr\{{\mathcal S}_V\}
= \Pr\li\{(1-\vep)\sigma_V^2 <
\overline{X_{2,N}^2} <
(1+\vep)\sigma_V^2\ri\} > 1-\de
\]
It follows that both lower and upper bounds in
(\ref{ineq2}) are true by taking $P_{s2}>0$ sufficiently
close to 0. Finally the proofs for
the rest of the theorem are omitted because
they are similar to those for the lower bounds
in (\ref{ineq1}) and (\ref{ineq2}).  \pfbox

\vspace{3mm}

\noindent {\bf Remark 5} \ Theorem \ref{Th3} shows that the
conservativeness of the coverage probability can be controlled
within $\f{\de}{2}$. Since the experimental effort grows in the
order of $\ln(\f{2}{\de})$ for a fixed relative width of the
confidence interval, a gain of $\f{\de}{2}$ for the coverage
probability can be obtained at the cost of increasing the
experimental effort by
\[
100 \li (\f{ \ln\left(
\f{2}{\de \sh 2}\right) } { \ln\left(\f{2}{\de}\right)} - 1 \ri ) \%
= \f{100\ln(2)}{\ln\left(\f{2}{\de}\right)}\%
\]
It can be seen that this percentage is very
insignificant for small $\de$. Therefore, the conservativeness of
the proposed interval estimation can be made very insignificant by
choosing large $P_{s1}$ and small $P_{s2}$.

\vspace{3mm}

The hypothesis $N_{\rm m2}>N_{\rm m1}$ in Theorem \ref{Th3} is
generically true for small $\delta$ and $\vep$ in light of the fact
that $N_{\rm m2}\ln\left(\f{2}{\de}\right) \approx N_{\rm
m1}\ln\left(\f{4}{\de}\right)$. Because the limiting case
$P_{s1}\rightarrow\infty$ is equivalent to the noiseless case, there
hold
\[
\lim_{P_{s1}\rightarrow\infty} \Pr\li\{\f{\overline{X_{1,N}^2}}
{(1+\vep)P_{s1}}<\sigma_H^2<\f{\overline{X_{1,N}^2}}
{(1-\vep)P_{s1}}\ri\} > 1-\de, \; \; \; \; \;
\lim_{P_{s1}\rightarrow\infty\atop{ P_{s2}\rightarrow 0}}
\Pr\li\{\f{\overline{X_{2,N}^2}}
{1+\vep}<\sigma_V^2<\f{\overline{X_{2,N}^2}}
{1-\vep}\ri\} > 1-\de
\]
provided that $N> \frac{2(3+2\vep)
\ln(\f{2}{\de})}{3\vep^2}\geq N_{\rm m1}$
where $\vep$ and $\de$ are positive numbers less than 1.
The sample complexity bound for the limiting
case improves the one in Theorem \ref{Th3}.

\subsection{Point Estimate}

In light of the results in \cite{per}, $\overline{X_{\ell,N}^2}$
($\ell=1$ or $\ell=2$) is the ML estimator for $P_{s\ell}\sigma_H^2
+ \sigma_V^2$ that is both sufficient and complete statistic. Let
$\wh{\si}_{H,N}^2$ and $\wh{\si}_{V,N}^2$ be point estimates for
$\sigma_H^2$ and $\sigma_V^2$, respectively based on
$\overline{X_{1,N}^2}$ and $\overline{X_{2,N}^2}$. Then
\be\label{tg} P_{s1}\wh{\si}_{H,N}^2 + \wh{\si}_{V,N}^2 =
\overline{X_{1,N}^2}, \; \; \; \; \; P_{s2}\wh{\si}_{H,N}^2 +
\wh{\si}_{V,N}^2 = \overline{X_{2,N}^2} \ee giving rise to the
following expressions for the point estimates (recall $\Delta P_s =
P_{s1}-P_{s2}>0$): \be\label{est} \wh{\si}_{H,N}^2 = \f{1}{\Delta
P_s} \li(\overline{X_{1,N}^2} - \overline{X_{2,N}^2}\ri), \; \; \;
\; \; \wh{\si}_{V,N}^2 =  \f{1}{\Delta P_s}
\li(P_{s1}\overline{X_{2,N}^2} - P_{s2}\overline{X_{1,N}^2}\ri) \ee
Our next result shows that both $\wh{\si}_{H,N}^2$ and
$\wh{\si}_{V,N}^2$ are optimum estimators for $\sigma_H^2$ and
$\sigma_V^2$, respectively.

\beT \la{Th5} Let $\{X_{1,k}\}_{k=1}^N$ and
$\{X_{2,k}\}_{k=1}^N$ be $2N$
{\rm i.i.d.} noisy samples with $P_{s1}>P_{s2}\geq 0$. Then
$\wh{\si}_{H,N}^2$ and $\wh{\si}_{V,N}^2$ as in $(\ref{est})$
are unbiased {\rm ML} and {\rm MV} estimates for $\sigma_H^2$ and $\sigma_V^2$,
respectively that achieve the {\rm Cram\'{e}r-Rao} lower bound.
Furthermore $\overline{X_{1,N}^2}$
and $\overline{X_{2,N}^2}$ as in $(\ref{tmp0})$ are
sufficient and complete statistics.
\eeT

\noindent Proof. It can be easily verified that
$\bb{E}[\overline{X_{\ell,N}^2}] = P_{s\ell}\sigma_H^2+\sigma_V^2$
for $\ell=1, 2$ implying that $\bb{E}[\wh{\si}_{H,N}^2] =
\sigma_H^2$ and $\bb{E}[\wh{\si}_{V,N}^2] = \sigma_V^2$.
Thus the estimator in $(\ref{est})$
is unbiased. By independence the joint PDF for $\{X_{1,k}\}_{k=1}^N$
and $\{X_{2,k}\}_{k=1}^N$ is given by \be\label{fxpdf}
f_X\left(\{X_{1,k}\}_{k=1}^N, \{X_{2,k}\}_{k=1}^N\right) =
f_X\left(\{X_{1,k}\}_{k=1}^N\right)f_X\left(\{X_{2,k}\}_{k=1}^N\right)
\ee by an abuse of notation where for either $\ell=1$ or $\ell=2$,
\[
f_X\left(\{X_{\ell,k}\}_{k=1}^N\ri) =
\f{1}{(P_{s\ell}\sigma_H^2+\sigma_V^2)^N}
{\rm exp}\left\{\sum_{i=1}^N \ln(|X_{\ell,i}|)\right\}
{\rm exp}\left\{-\f{N\overline{X_{\ell,N}^2}}
{P_{s\ell}\sigma_H^2+\sigma_V^2}\right\}
\]
Taking natural logarithm on both sides yields
\[
\ln(f_X) = -\ln(P_{s1}\sigma_H^2+\sigma_V^2)^N
- \ln(P_{s2}\sigma_H^2+\sigma_V^2)^N + \sum_{i=1}^N
\ln(|X_{1,i}X_{2,i}|)
- \f{N\overline{X_{1,N}^2}}
{P_{s1}\sigma_H^2+\sigma_V^2}-\f{N\overline{X_{2,N}^2}}
{P_{s2}\sigma_H^2+\sigma_V^2}
\]
Setting partial derivatives of $\ln(f_X)$ with respect to
$\sigma_H^2$ and $\sigma_V^2$ to zeros lead to the
following two equations:
\bee
\f{P_{s1}\overline{X_{1,N}^2}}{(P_{s1}\sigma_H^2+\sigma_V^2)^2} +
\f{P_{s2}\overline{X_{2,N}^2}}{(P_{s2}\sigma_H^2+\sigma_V^2)^2} &=&
\f{P_{s1}}{P_{s1}\sigma_H^2+\sigma_V^2} +
\f{P_{s2}}{P_{s2}\sigma_H^2+\sigma_V^2}\\
\f{\overline{X_{1,N}^2}}{(P_{s1}\sigma_H^2+\sigma_V^2)^2} +
\f{\overline{X_{2,N}^2}}{(P_{s2}\sigma_H^2+\sigma_V^2)^2} &=&
\f{1}{P_{s1}\sigma_H^2+\sigma_V^2} +
\f{1}{P_{s2}\sigma_H^2+\sigma_V^2}
\eee
Solving for $\overline{X_{1,N}^2}$ and $\overline{X_{2,N}^2}$
from the above two equations coincide with the
expressions in (\ref{tg}) at
$\sigma_H^2 = \wh{\sigma}_{H,N}^2$ and $\sigma_V^2
= \wh{\sigma}_{V,N}^2$ that in turn gives the
point estimates in (\ref{est}).
The uniqueness of the solution implies that the
points estimates in (\ref{est}) are indeed ML
that maximize the PDF $f_X(\cdot)$ for every pair of sets of
$N$ i.i.d. noisy measurement samples satisfying the hypothesis.
Denote $\theta = \left[\begin{array}{cc}
\sigma_H^2 & \sigma_V^2\end{array}\right]^T$
as the parameter vector for estimation where
superscript $^T$ denotes transpose. Then the
corresponding FIM (Fisher information matrix)
can be shown to be, after lengthy calculation,
\[
{\rm Fim}(\sigma_H^2, \sigma_V^2)
:= \bb{E}\left\{\left[\f{\partial\ln(f_X)}{\partial\theta}\right]
\left[\f{\partial\ln(f_X)}{\partial\theta}\right]^T\right\} =
N\left[\begin{array}{cc}
\f{P_{s1}^2}{\sigma_1^4} + \f{P_{s2}^2}{\sigma_2^4} &
\f{P_{s1}}{\sigma_1^4} + \f{P_{s2}}{\sigma_2^4} \\
\f{P_{s1}}{\sigma_1^4} + \f{P_{s2}}{\sigma_2^4} &
\f{1}{\sigma_1^4} + \f{1}{\sigma_2^4}
\end{array}\right]
\]
where for simplicity in notations, $\sigma_{\ell}
= P_{s\ell}\sigma_H^2+\sigma_V^2$ is used with
$\ell=1, 2$. Recall that the inverse
of ${\rm Fim}(\sigma_H^2, \sigma_V^2)$
is the smallest error covariance achievable by all
unbiased estimators. In arriving to the above expression, the
fact that $\bb{E}[G^4]=3\sigma^4$ for $G$ distributed as
$\mscr{N}(0, \sigma)$, and the observations:
\bee
\bb{E}\li[\overline{X_{\ell,N}^2}\ri]
&=& P_{s\ell}\sigma_H^2+\sigma_V^2, \; \; \; \; \;
\bb{E}\li[(\overline{X_{1,N}^2})(\overline{X_{2,N}^2})\ri]
= (P_{s1}\sigma_H^2+\sigma_V^2)(P_{s2}\sigma_H^2+\sigma_V^2)\\
\bb{E}\li[\left(\overline{X_{\ell,N}^2}\right)^2\ri]
&=& (P_{s\ell}\sigma_H^2+\sigma_V^2)^2 +
\f{1}{N}(P_{s\ell}\sigma_H^2+\sigma_V^2)^2,
\; \; \; \; \; \ell=1, 2
\eee
are crucial. Now with the point estimates in (\ref{est}),
the corresponding error covariance can be shown to be
\bee
{\rm Cov}(\sigma_H^2, \sigma_V^2) &=& \bb{E}
\li\{\left[\begin{array}{c}
\wh{\si}_{H,N}^2 - \sigma_H^2 \\
\wh{\si}_{V,N}^2 - \sigma_V^2
\end{array}\right]\left[\begin{array}{cc}
\wh{\si}_{H,N}^2 - \sigma_H^2 &
\wh{\si}_{V,N}^2 - \sigma_V^2
\end{array}\right]\ri\} \\
&=& \f{1}{N\Delta P_s}
\left[\begin{array}{cc}
\sigma_1^4+\sigma_2^4 & -(P_{s1}\sigma_2^4+P_{s2}\sigma_1^4)\\
-(P_{s1}\sigma_2^4+P_{s2}\sigma_1^4) &
P_{s1}^2\sigma_2^4+P_{s2}^2\sigma_1^4
\end{array}\right]
\eee
where again the notation $\sigma_{\ell}^2 = P_{s\ell}\sigma_H^2+\sigma_V^2$
with $\ell=1, 2$ is used. Consequently there holds the identity
${\rm Fim}(\sigma_H^2, \sigma_V^2){\rm Cov}(\sigma_H^2, \sigma_V^2)
= I$, validating the fact that the Cram\'{e}r-Rao lower bound
is indeed achieved by the point estimates in (\ref{est}) that are indeed
MV estimates. Finally the
sufficient and complete statistic of $\overline{X_{1,N}^2}$
and $\overline{X_{2,N}^2}$ can be shown in the same way as in \cite{per}.
Specifically the joint PDF in (\ref{fxpdf}) has the form
\[
f_X = g\left(\overline{X_{1,N}^2}, \overline{X_{2,N}^2},
\sigma_H^2, \sigma_V^2\right)h\left(\{X_{1,k}\}_{k=1}^N,
\{X_{2,k}\}_{k=1}^N\right),
\; \; \; \; \; h(\cdot, \cdot) =
\prod_{i=1}^N |X_{1,i}X_{2,i}|
\]
Hence in light of the Neyman factorization criterion,
$\overline{X_{1,N}^2}$ and $\overline{X_{2,N}^2}$ are
sufficient statistics. The proof for the complete statistics
is again the same as that in \cite{per} that is skipped. \pfbox

\vspace{3mm}

The result on the Cram\'{e}r-Rao lower bound for finite
$N$ in Theorem \ref{Th5} is surprising
which is normally achievable only asymptotically. Next we
present several results for the
point
estimates in (\ref{est}) regarding the sample complexity bound.
We begin with a corollary that follows from the
results on interval estimates in Subsection 4.1.

\begin{corollary} \la{Th4} Under the same hypotheses as in {\rm Theorem}
$\ref{Th3}$, there holds
\be\label{theq1} 
1-\de\leq  \Pr \li \{  \li | \f{ \wh{\si}_{H,N}^2  - \si_H^2 } {
\si_H^2  } \ri | < \vep \ri \} < 1 - \f{\de}{2} \ee for sufficiently
large $P_{s1}$ with constant $P_{s2}$. In addition \be\label{theq2}
1-\de\leq  
\Pr \li \{  \li | \f{
\wh{\si}_{V,N}^2  - \si_V^2 }
{ \si_V^2  } \ri | < \vep \ri \} < 1 - \f{\de}{2}
\ee
for sufficiently small $P_{s2}$ with constant $P_{s1}$.
\end{corollary}

\noindent Proof. We prove only (\ref{theq1}) as the
proof for (\ref{theq2}) is similar that will be omitted.
Define the set of events:
\[
\wh{\mathcal S}_H := \li\{\sigma_H^2: \; \;
\frac{1}{\Delta P_s(1+\vep)}\li(\overline{X_{1,N}^2}
- \overline{X_{2,N}^2}\ri)
<\sigma_H^2< \frac{1}{\Delta P_s(1-\vep)}\li(\overline{X_{1,N}^2}
- \overline{X_{2,N}^2}\ri)\ri\}
\]
Then $\wh{\mathcal S}_H\subseteq {\mathcal S}_H$, and thus
$\Pr\{\wh{\mathcal S}_H\} \leq \Pr\{{\mathcal S}_H\}$.
More importantly $\wh{\mathcal S}_H \rightarrow {\mathcal S}_H$
as $P_{s1}\rightarrow \infty$. Therefore
\[
\Pr\{\wh{\mathcal S}_H\} = \Pr\li\{
\frac{\wh{\si}_{H,N}^2}{1+\vep}
<\sigma_H^2< \frac{\wh{\si}_{H,N}^2}{1-\vep}\ri\}\\
=  \Pr \li \{  \li | \f{
\wh{\si}_{H,N}^2  - \si_H^2
} { \si_H^2  } \ri | < \vep \ri \}
\ \rightarrow \ \Pr\{{\mathcal S}_H\}
\]
as $P_{s1}\rightarrow \infty$.
A similar argument as in the proof of Theorem \ref{Th3}
can thus be used to conclude (\ref{theq1}). \pfbox

\vspace{3mm}

Similar to the discussion after Remark 5, we have that
for $N >\f{ 2(3 + 2 \vep) \ln \f{2}{\de}   } { 3 \vep^2  }$, there hold
\be\label{tmp02}
\lim_{P_{s1}\to \iy} \Pr \li \{  \li | \f{ \wh{\sigma}_{H,N}^2 - \si_H^2
} { \si_H^2 } \ri | < \vep \ri \} > 1 - \de
\; \; \; \; \;
\lim_{P_{s2} \to 0} \Pr \li \{  \li | \f{ \wh{\sigma}_{V,N}^2 - \si_V^2 }
{ \si_V^2  } \ri | < \vep \ri \} > 1 - \de
\ee
where $\vep$ and $\de$ are positive numbers less than
$1$.
The next result is concerned with the point estimate related to SNR:
\be
\wh{\si}_{\f{H}{V},N}^2 =  \f{\wh{\si}_{H,N}^2}{\wh{\si}_{V,N}^2}
= \f{\overline{X_{1,N}^2} - \overline{X_{2,N}^2}}
{P_{s1}\overline{X_{2,N}^2} - P_{s2}\overline{X_{1,N}^2}}
\ee

\beT \la{Th6} Let $\vep$ and $\de$ be positive numbers less than
$1$.  Let $P_{s2} = 0$.  Suppose that $N$ is the minimum
positive integer such that
$F_{\chi^2, 2N} \li ( 2N(1+ \f{\vep}{2 + \vep} ) \ri ) -
F_{\chi^2, 2N} \li ( 2N(1 - \f{\vep}{2 + \vep}) \ri ) \geq \sq{1 -
\de}$ where $F_{\chi^2, 2N} (\cdot)$ denotes the
{\rm CDF} of $\chi^2$ random variable with $2N$ degrees of freedom. Then
\[
\lim_{P_{s1} \to \iy} \Pr\li \{
\li| \f{\wh{\sigma}_{\f{H}{V},N}^2 - \f{\si_H^2}{\si_V^2}}
{\f{\si_H^2}{\si_V^2}} \ri| < \vep
\ri \} > 1 - \de
\]
\eeT

\noindent Proof. Define two normalized random variables by
\be
\wh{\theta}_{H,N} := \frac{\overline{X_{1,N}^2}}
{P_{s1}\sigma_H^2+ \sigma_V^2},
\; \; \; \; \;
\wh{\theta}_{V,N} := \frac{\overline{X_{2,N}^2}}{\sigma_V^2}
\ee
Straightforward algebraic manipulations give the following
chain of equations:
\bee
\Pr\li \{
\li| \f{\wh{\sigma}_{\f{H}{V},N}^2 - \f{\si_H^2}{\si_V^2}}
{\f{\si_H^2}{\si_V^2}} \ri| < \vep \ri \} &=&
\Pr\li \{\f{\wh{\sigma}_{\f{H}{V},N}^2}{1+\vep}
< \f{\si_H^2}{\si_V^2} <
\f{\wh{\sigma}_{\f{H}{V},N}^2}{1-\vep}\ri \} \\
&=& \Pr\li \{\f{1}{P_{s1}(1+\vep)}\left(
\f{\overline{X_{1,N}^2}}{\overline{X_{2,N}^2}}-1\right)
< \f{\si_H^2}{\si_V^2} < \f{1}{P_{s1}(1-\vep)}
\left(\f{\overline{X_{1,N}^2}}{\overline{X_{2,N}^2}}-1\right)\ri\}\\
&=& \Pr\li \{\left(\f{P_{s1}(1-\vep)\sigma_H^2}{\sigma_V^2}+1\right)
< \f{\overline{X_{1,N}^2}}{\overline{X_{2,N}^2}} <
\left(\f{P_{s1}(1+\vep)\sigma_H^2}{\sigma_V^2}+1\right)\ri\} \\
&=& \Pr\li\{\f{P_{s1}(1-\vep)\sigma_H^2+\sigma_V^2}
{P_{s1}\sigma_H^2+\sigma_V^2}
< \f{\wh{\theta}_{H,N}}{\wh{\theta}_{V,N}} <
\f{P_{s1}(1+\vep)\sigma_H^2+\sigma_V^2}
{P_{s1}\sigma_H^2+\sigma_V^2}\ri\} \\
&\rightarrow& \Pr\li\{1-\vep
< \f{\wh{\theta}_{H,N}}{\wh{\theta}_{V,N}} <
1+\vep\ri\}
\eee
as $P_{s1}\rightarrow \infty$.
To complete the proof, define two sets of events as follows:
\bee
{\mathcal S}_{\theta} &:=& \li\{(\wh{\theta}_{H,N},
\wh{\theta}_{V,N}): \; \; \;
(1-\vep)\wh{\theta}_{V,N}
< \wh{\theta}_{H,N} < (1+\vep)\wh{\theta}_{V,N}\ri\},
\; \; \; \; \; g(\vep) := \frac{\vep}{2+\vep} \\
\wh{\mathcal S}_{\theta} &:=&  \li\{(\wh{\theta}_{H,N},
\wh{\theta}_{V,N}): \; \; \;
1-g(\vep) < \wh{\theta}_{H,N} < 1+g(\vep) \; \; \; \& \; \; \;
1-g(\vep) < \wh{\theta}_{V,N} < 1+g(\vep)\ri\}
\eee
It can be shown that $\wh{\mathcal S}_{\theta}
\subseteq {\mathcal S}_{\theta}$. Indeed in
the plane of $(\wh{\theta}_{H,N}, \wh{\theta}_{V,N})$,
the set $\wh{\mathcal S}_{\theta}$ defines a square
area inside the sector area defined by the set
${\mathcal S}_{\theta}$. More specifically
\[
\wh{\mathcal S}_{\theta} \subseteq
\li\{(\wh{\theta}_{H,N}, \wh{\theta}_{V,N}):
\; \; \;
1-\f{\vep}{2-\vep} <
\wh{\theta}_{H,N} <1+g(\vep)
\; \; \; \& \; \; \; 1-g(\vep) < \wh{\theta}_{V,N} < 1+\f{\vep}{2-\vep}\ri\}
\subseteq {\mathcal S}_{\theta}
\]
It follows that as $P_{s1}\rightarrow \infty$, there holds
\bee
\Pr\{{\mathcal S}_{\theta}\} &=&
\Pr\li\{1-\vep
< \f{\wh{\theta}_{H,N}}{\wh{\theta}_{V,N}} <
1+\vep\ri\} \geq \Pr\{\wh{\mathcal S}_{\theta}\} \\
&=& \Pr\li\{1-g(\vep) < \wh{\theta}_{H,N} < 1+g(\vep)\ri\}
\Pr\li\{1-g(\vep) < \wh{\theta}_{V,N} < 1+g(\vep)\ri\}>1-\de
\eee
by independence of the events of $\wh{\theta}_{H,N}$
and $\wh{\theta}_{H,N}$, and by the hypothesis on $N$.
\pfbox

\vspace{3mm}

By noting that $F_{\chi^2, 2N}
\li (2N(1+ \f{\vep}{2 + \vep} ) \ri ) -
F_{\chi^2, 2N} \li (2N(1 - \f{\vep}{2 + \vep}) \ri ) \geq \sq{1 -
\de}$, if $N > \f{ 2(2 + \vep)(6+5\vep)\ln
\f{4}{\de}} { 3 \vep^2 }$, the next result follows from
Theorem \ref{Th6}.

\begin{corollary}\la{Th7} Let $\vep$ and $\de$ be positive numbers less than
$1$. Suppose that $P_{s2} = 0$.
Then
\[
\lim_{P_{s1} \to \iy} \Pr\li \{
\li| \f{\wh{\sigma}_{\f{H}{V},N}^2 - \f{\si_H^2}{\si_V^2}}
{\f{\si_H^2}{\si_V^2}} \ri| < \vep
\ri \} > 1 - \de
\]
provided that $N
> \f{ 2(2 + \vep)(6+5\vep)\ln
\f{4}{\de}} { 3 \vep^2 } \approx \f{8\ln(\f{4}{\de})}{\vep^2}$ for small $\vep$.
\end{corollary}

We note that $\f{ 2(2 + \vep)(6+5\vep)\ln
\f{4}{\de}} { 3 \vep^2 } \approx \f{8\ln(\f{4}{\de})}{\vep^2}$ for small $\vep$.
Thus the {\it a priori} bounds for the sample
complexity in Theorem \ref{Th6} and
Corollary \ref{Th7} are roughly 4 times to that in Corollary
\ref{Th4}. It signifies the difficulty in estimation of SNR
in Rayleigh channels, and indicates that
the margin of error manifested by $\vep$ is more expensive
to reduce as compared with the level of confidence
manifested by $\de$.

\section{Numerical Simulations}

This section presents numerical results to illustrate
our proposed statistical theory in measurements and
estimation of Rayleigh fading channels. For simplicity
$\sigma_H^2=1$ and $P_{s1} = 1$ are assumed throughout the
section that can always be
made true by a suitable normalization. Hence
the SNR is the same as $\f{\sigma_H^2}{\sigma_V^2}$,
and large $P_{s1}$ reduces to large SNR. Moreover
we consider only the case $P_{s2}=0$.

\vspace{1mm}

Because the noiseless case has been shown in Figure 1,
we begin with the noisy measurement samples
under SNR $=$ 20 dB by generating
$\{X_{1,i}\}_{i=1}^N$ and $\{X_{2,i}\}_{i=1}^N$
for different sample size $N$. The simple
averages $\overline{X_{\ell,N}^2}$ as in (\ref{tmp0})
are then calculated for $\ell=1, 2$. We choose
$\de=0.01$ for the associated confidence level. In
light of Theorem \ref{Th2}, there holds the joint probability
\[
\Pr\li\{A_1<\sigma_H^2<A_2 \; \& \;
B_1<\sigma_V^2<B_2\ri\} > 1-\de
\]
where $A_1, A_2, B_1, B_2$ are functions of
$\vep$, $N$, and $\overline{X_{\ell,N}^2}$ for $\ell=1, 2$
which represent the interval estimates.
In Figure 3 we plotted the ratios
\[
\frac{A_2}{A_1} = \f{(1+\vep)\overline{X_{1,N}^2}
- (1-\vep)\overline{X_{2,N}^2}}
{(1-\vep)\overline{X_{1,N}^2}
- (1+\vep)\overline{X_{2,N}^2}}, \; \; \; \; \; \; \;
\frac{B_2}{B_1} = \f{P_{s1}(1+\vep)\overline{X_{2,N}^2}
- P_{s2}(1-\vep)\overline{X_{1,N}^2}}
{P_{s1}(1-\vep)\overline{X_{2,N}^2}
- P_{s2}(1+\vep)\overline{X_{1,N}^2}}
\]
for the case $P_{s2}=0$ with SNR $=$ 20 dB (left) and SNR $=$ 60 dB
(right). It is interesting to observe that although
$\overline{X_{\ell,N}^2}$, $\ell=1, 2$, are random for each $N$, the
above ratios are almost deterministic owing to the small $\de$ value
used. In addition the SNR values affect little for the two ratios at
large $N$ values. We also plotted the noiseless ratio in dashed line
as a comparison that should serve as a lower limit.

\vspace*{2.5in}
\begin{figure}[ht]\label{fig4}
\includegraphics{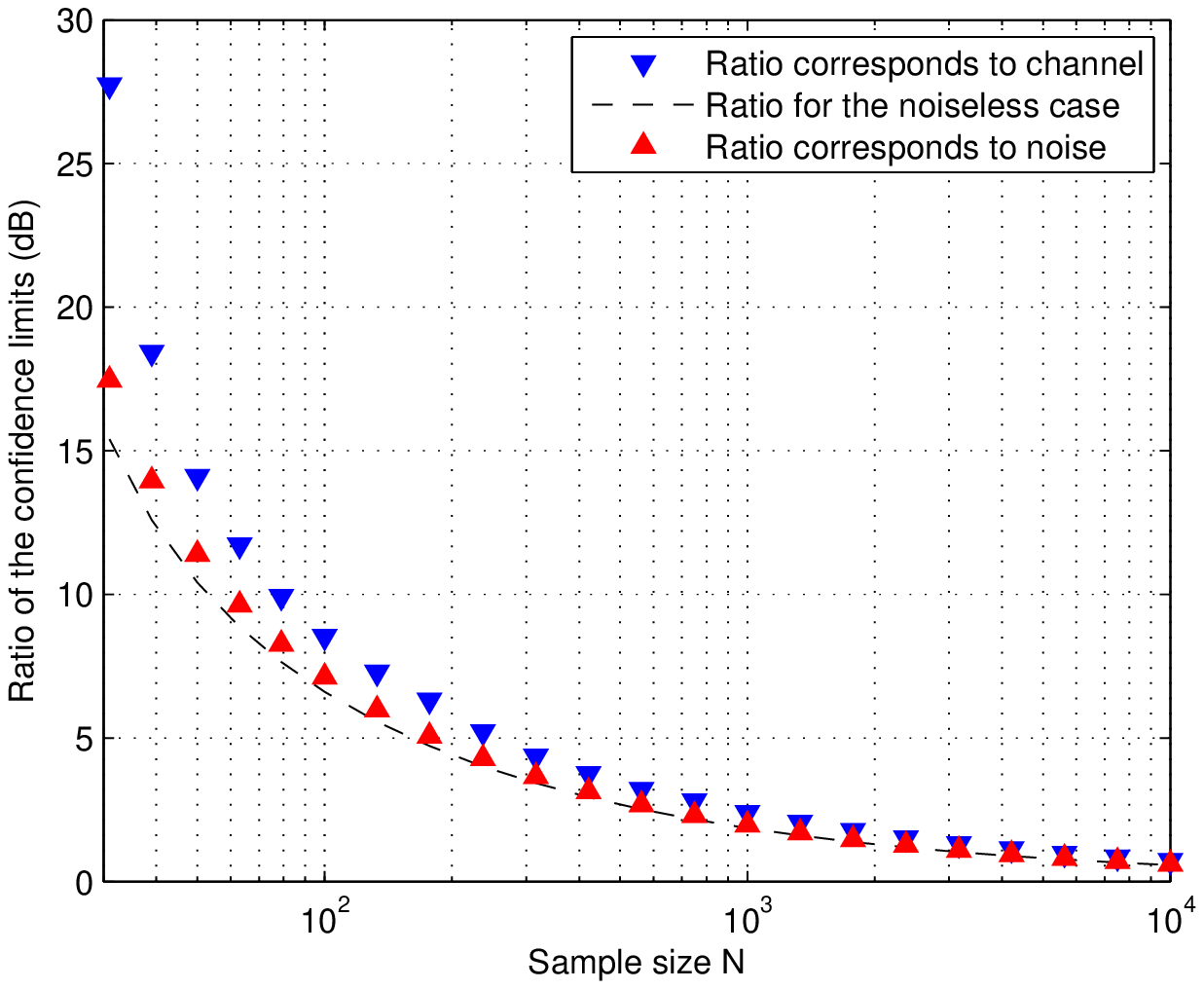}
\includegraphics{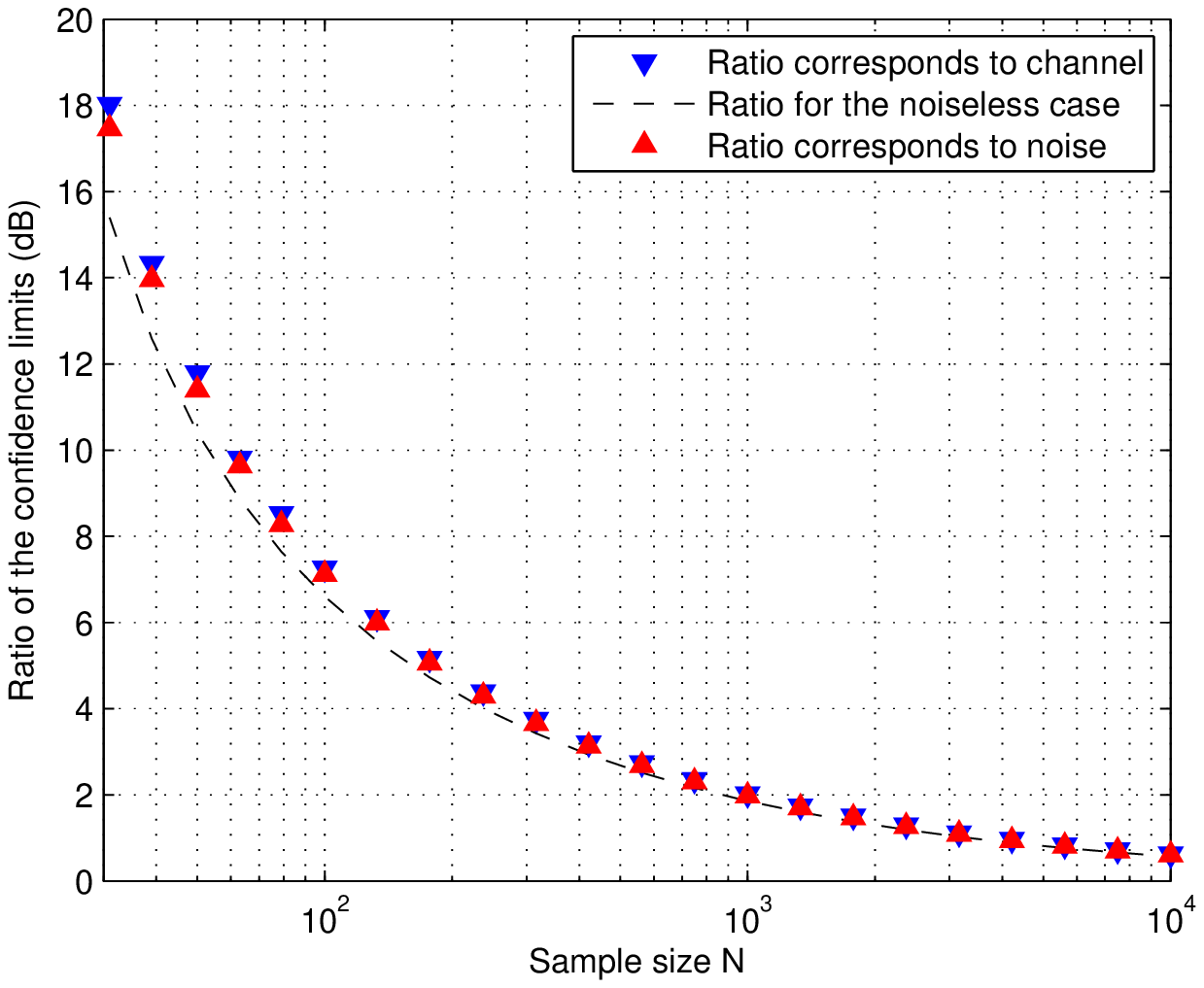}
\caption{Ratios of the upper/lower confidence limits
with SNR $=$ 10dB (left) and SNR $=$ 20dB (right)}
\end{figure}

For estimation of the SNR, there holds $\Pr\{D_1<\f{\sigma_H^2}{\sigma_V^2}<D_2\}$
in light of Corollary \ref{co2} with
\[
\f{D_2}{D_1} = \f{(1+\vep)\overline{X_{1,N}^2} -
(1-\vep)\overline{X_{2,N}^2}}{(1-\vep)\overline{X_{1,N}^2} -
(1+\vep)\overline{X_{2,N}^2}}
\]
In Figure 4, the above is plotted against the sample size $N$
for the cases SNR = 10 and SNR $=$ 60.

\vspace*{2.2in}
\begin{figure}[ht]\label{fig3}
\includegraphics{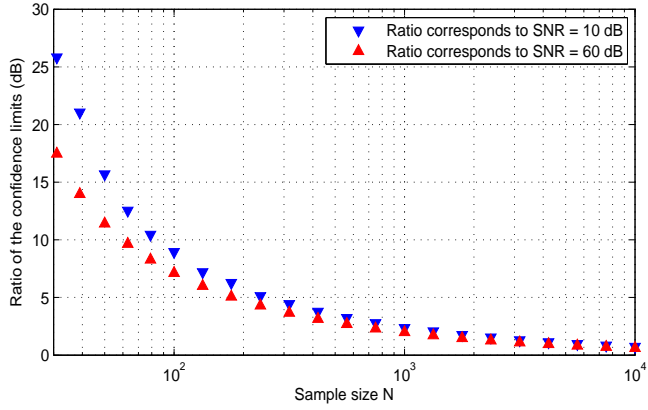}
\caption{Ratios of the upper/lower confidence limits
corresponds with $\de=0.01$}
\end{figure}

It is commented that for both Figure 3 and Figure 4, the
{\it a priori} bound $N>\f{2(3+2\vep)\ln(\f{4}{\de})}{3\vep^2}$,
rather than the relation in (\ref{nmin2}), is used.
The results are nevertheless close to each other. In Figure 5,
we plotted $\vep$ vs. $N$ (left) for the case $\de = 0.05$,
and plotted $\de$ vs. $N$ (right) for the $\vep = 0.05$.
In both cases, {\it a priori} bounds are used
in both the noisy and noiseless cases,
plus the use of $\chi^2$ in the noiseless case
that is governed by the relation in (\ref{nmin1}).

\vspace*{2.4in}
\begin{figure}[ht]\label{fig5}
\includegraphics{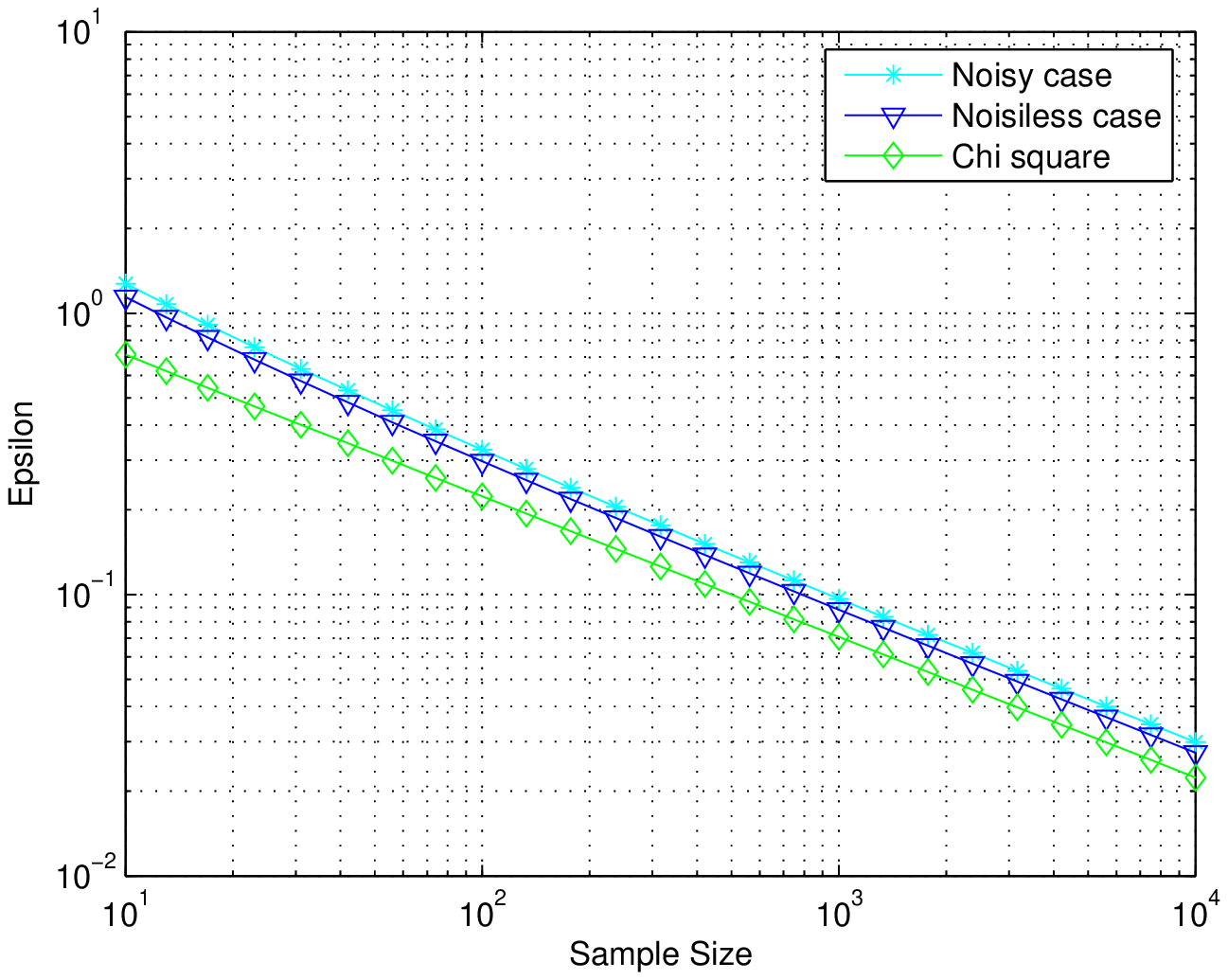}
\includegraphics{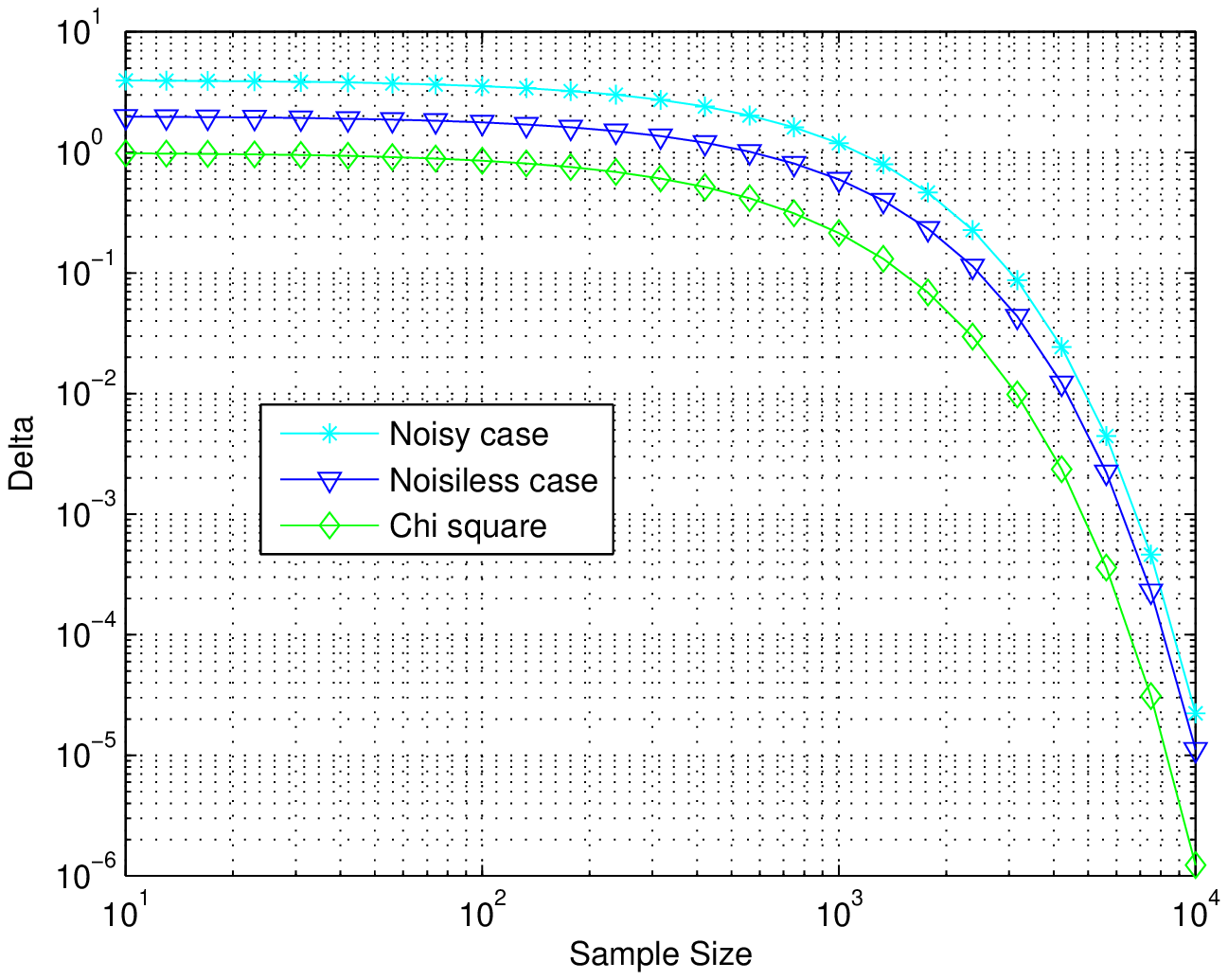}
\caption{Plots of $\vep$ vs. $N$  with $\de=0.05$
(left),  and $\de$ vs. $N$ with $\vep=0.05$ (right)}
\end{figure}

As expected, the curve based on $\chi^2$ CDF in the noiseless case
serves as a lower bound. In addition it is observed that the two
curves based on the {\it a priori} bounds are almost identical for
the one on left, and close to each other for the one on right. This fact
indicates that our results on measurement and estimation of noisy
Rayleigh fading channels are not conservative. It is further
observed that the $\vep$ curves decrease at a constant slope in
log-log scale, and the slope of decrease is rather slow with respect
to the sample size $N$. This fact indicates that the reduction of
$\vep$ (with fixed $\de$) is expensive in terms of increasing $N$. On
the other hand the $\de$ curves decrease at accelerated slopes with
respect to the sample size $N$, implying that the reduction of $\de$
(with fixed $\vep$) is relatively cheap, especially at large sample
size $N$.
Figure 5 validates that
the measurement sample size is roughly inversely
proportional to the square of the margin of
error $\vep$ and is
linear with respect to the logarithm
of the inverse of the gap $\delta$.

\vspace{1mm}

The last simulation example is worked out to demonstrate the optimum
estimator obtained in Theorem \ref{Th5}. Different from the previous
two cases, $M (>>1)$ sets of measurement samples are taken to assess
the average performance for the underlying statistical estimation.
For each sample size $N$, $\{X_{\ell,i}(k)\}_{i=1}^N$ is generated
for $\ell=1, 2$, and $1\leq k\leq M$ with $M=500$. The optimum
estimator in (\ref{est}) is used to compute estimates
$\wh{\sigma}_{H,N}^2(k)$ and $\wh{\sigma}_{V,N}^2(k)$. The
estimation error is then averaged to yield
\[
{\rm RMES}_H \approx \sqrt{\f{1}{M}\sum_{k=1}^M
\li|\wh{\sigma}_{H,N}^2(k)-\sigma_H^2\ri|^2}, \; \; \; \; \;
{\rm RMES}_V \approx \sqrt{\f{1}{M}\sum_{k=1}^M
\li|\wh{\sigma}_{V,N}^2(k)-\sigma_V^2\ri|^2}
\]
that are plotted together with the Cram\'{e}r-Rao lower bounds:
\[
{\rm CRB}_H = \sqrt{\f{(P_{s1}\sigma_H^2+\sigma_V^2)^2 +
(P_{s2}\sigma_H^2+\sigma_V^2)^2}{N\Delta P_s}},
\; \; \; \; \; {\rm CRB}_V = \sqrt{\f{P_{s2}^2(P_{s1}\sigma_H^2+\sigma_V^2)^2 +
P_{s1}^2(P_{s2}\sigma_H^2+\sigma_V^2)^2}{N\Delta P_s}}
\]
versus the sample size $N$ as shown in Figure 6 where $P_{s2}=0$ is
used. It can be seen that the RMSEs coincide very well with the
CRBs, validating the optimality of the estimator in (\ref{est}).
Again the change of SNR affects little on estimation of $\sigma_H^2$
but changes the RMSE for $\sigma_V^2$ that is mainly due to the
change of $\sigma_V^2$ by a factor of 10. In the case $M<500$, the
RMSE lines are less straight and fluctuate more as $M$ becomes
smaller, but the overall trend holds.

\vspace*{2.6in}
\begin{figure}[ht]\label{fig2}
\includegraphics{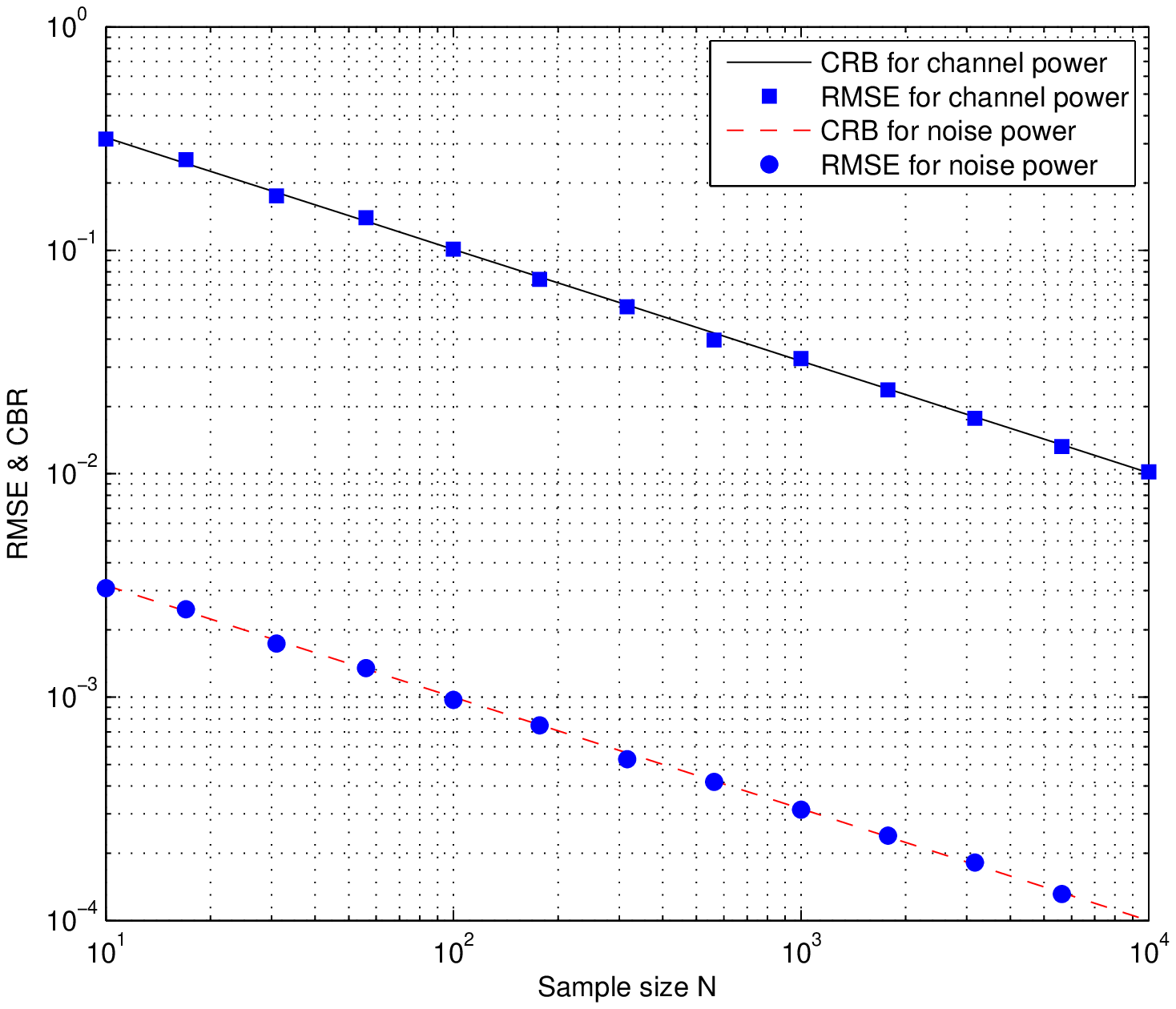} \includegraphics{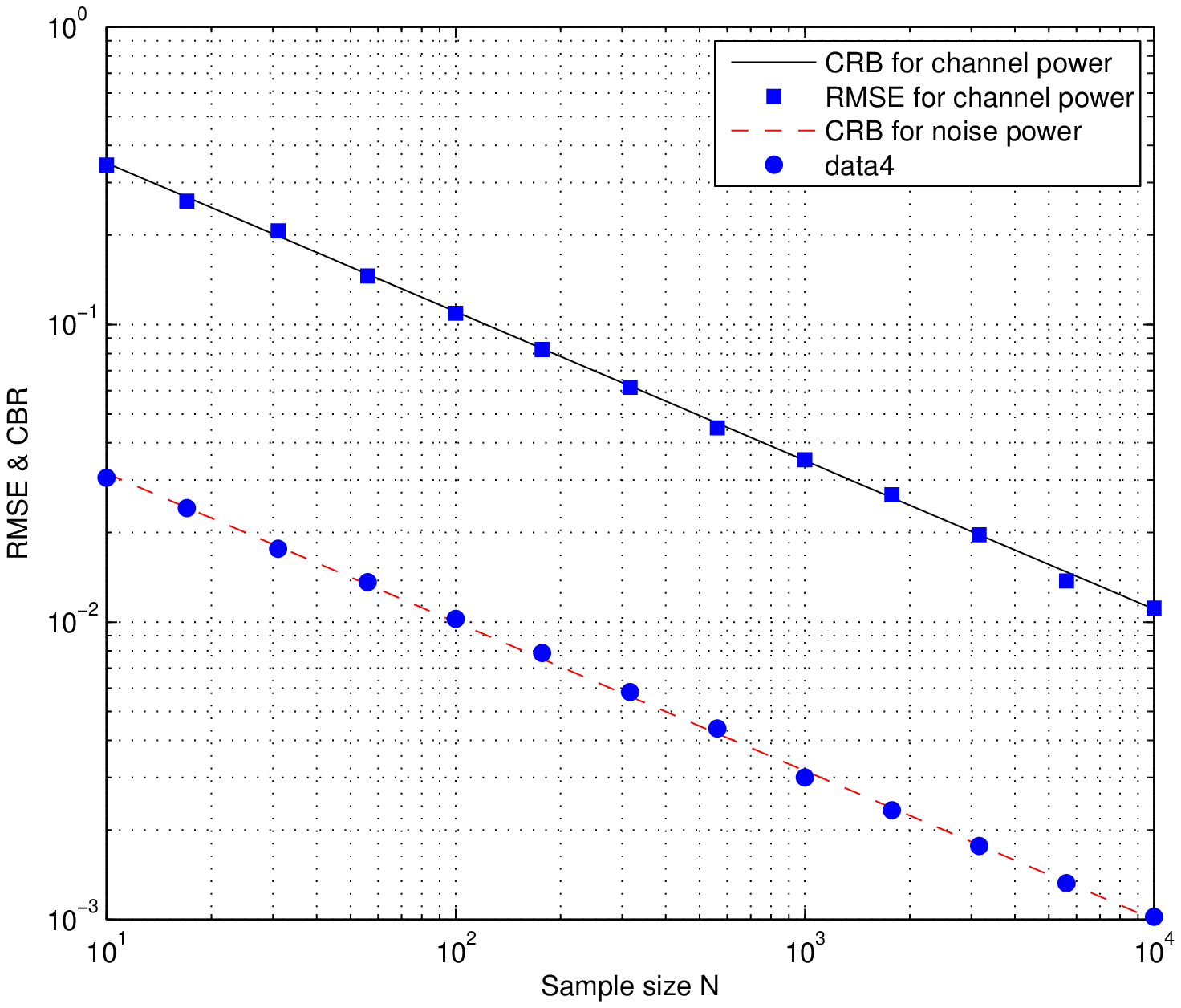} \caption{RMSE and CRB vs. $N$ with SNR $=$
20dB (left) and SNR $=$ 10dB (right)}
\end{figure}

\section{Conclusion}

Statistical estimation of Rayleigh fading channels has been
investigated based on both noiseless and noisy measurement samples.
Complete solutions are derived for the associated sample complexity
problem and provided for the optimum estimator problem in
measurement and estimation of Rayleigh fading channels. Specifically
our {\it a priori} bounds on measurement sample sizes ensure the
prescribed margin of error and confidence level and are contrast to
the existing work reported in the literature. In dealing with the noisy
measurement samples, our proposed novel signaling scheme with two
different signal strengths is instrumental in extracting the
statistical information on mean channel power, noise power, and SNR. Such
a novel signaling scheme enables us to derive the sample complexity
bounds for both interval and point estimates that are tight for the
mean
channel power and noise power albeit less tight for the SNR. More
importantly it leads to the optimum estimator that is both an ML and
MV estimator and that achieves the Cram\'{e}r-Rao lower bound. The
results presented in this paper constitute an independent
statistical theory for measurement and estimation of Rayleigh fading
channels. It should be emphasized that the sample
complexity solution in the noiseless case is also instrumental
without which the results for the case of noisy measurements
are not possible. The numerical simulations illustrate that our proposed
statistical theory is effective in statistical
estimation of Rayleigh fading channels. Specifically the
simulation examples indicate that our results based on
the noisy measurement samples
are close to that based on the
noiseless measurement samples, and the signal power is not
required to be high that can be compensated for by using large sample
size. Currently we are investigating the Nakagami fading channel
which in a special case reduces to the Rayleigh fading channel,
and aim at extension of our results
on the measurement sample size and optimum
estimator. It
is our objective to apply our results
to other more general
wireless fading channels and to broad
applications of our proposed
statistical theory.


\begin{thebibliography}{10}

\bibitem{ak} A. Abdi and M. Kaveh, ``Performance comparison
of three different estimators for the Nakagami $m$ parameter
using Monte Carlo simulation,'' {\it IEEE Commun. Lett.},
vol. 4, pp. 119-121, April 2000.

\bib{b63} P.A. Bello, ``Characterization of
randomly time-variant linear channels,'' {\it IEEE Trans. Commun.},
vol. 11, pp. 360--393, Dec. 1963.

\bib{besw} J. Beek, O. Edfors, M. Sandell, S. Wilson, and P. B\"{o}rjesson,
``On channel estimation in OFDM systems,'' in {\it Proceedings of
IEEE VTC 95}, Chicago, 1995.

\bibitem{b} M. Barbiloni, C. Carciofi, G. Falciasecca, M. Frallone,
and P. Grazioso, ``A measurement-based methodology for the
determination of validity domains of prediction models in urban
environments,'' {\it IEEE Trans. Veh. Technol.}, vol. 49, pp.
1508-1515, Sept. 2000.

\bib{chernoff} H. Chernoff, ``A measure of asymptotic efficiency for tests
of a hypothesis based on the sum of observations,'' {\it Annals of
Mathematical Statistics}, vol. 23, no. 4, pp. 493-507, 1952.

\bibitem{cb} J. Cheng and N.C. Beaulieu, ``Generalized
moment estimators for the Nakagami fading parameters,''
{\it IEEE Commun. Lett.}, vol. 6, pp. 144-146, April 2002.

\bib{chen} Y. Chen and N.C. Beaulieu, ``Estimators using noisy channel
samples for fading distribution parameters,'' {\it IEEE Trans.
Commun.}, vol. 53, no. 8, Aug. 2005, pp. 1274-1277.

\bibitem{c} N. Cotanis, ``Estimating radio coverage for new
mobile wireless services data collection and pre-processing,''
ICT2001, Bucharest, Romania, June 2001.

\bib{czink} N. Czink, G. Matz, D. Seethaler, and F.
Hlawatsch, ``Improved MMSE estimation of correlated MIMO channels
using a structured correlation estimator,'' {\it Proc. IEEE
SPAWC-2005}, New York, June 2005, pp. 595-599.

\bib{dietrich} F.A. Dietrich, T. Ivanov, and W.
Utschick, ``Estimation of channel and noise correlations for MIMO
channel estimation,'' {\it Proc. of the Workshop on Smart Antennas},
Germany, 2006.

\bib{dogandzic} A. Dogandzic and J. Jin, ``Estimating
statistical properties of MIMO fading channels,'' {\it IEEE
Transactions on Signal Processing}, vol. 53, pp. 3065-3080, Aug. 2005.

\bib{gghn} G. Gu, J. He, X. Gao, and M. Naraghi-Pour,
``An analytic approach to modeling and estimation of OFDM
channels,'' in {\it Proceedings of IEEE 2004 Global Communications
Conference} (Dallas, TX), Nov. 2004, and to appear in
{\it IEEE Trans. Veh. Technol.}, July 2007.

\bib{per} M. Peritsky, ``Statistical estimation of mean signal strength
in a Rayleighfading environment,'' {\it IEEE Trans. Commun.}, vol.
21, no. 11, Nov. 1973, pp. 1207-1213.

\bibitem{p} J.D. Parsons, {\it The mobile Radio Propagation
Channel}, Pentech House, 1994.

\bibitem{proakis} J.G. Proakis, {\it Digital Communications}, McGraw-Hill,
2000.

\bibitem{r} T.S. Rappaport, {\it Wireless Communications --- Principles
and Practice}, Prentice-Hall, Upper Saddle River, NJ, 1999.

\bibitem{sa} M.K. Simon and M.-S. Alouini, ``Average bit-error
probability performance for optimum diversity combining
of noncoherent FSK over Rayleigh fading channels,'' {\it IEEE
Trans. Commun.}, vol. 51, pp. 566-569, April 2003.

\bibitem{st} G.L. St$\ddot{\rm u}$ber, {\it Principles of
Mobile Communication}, second edition, Kluwer Academic Publishers,
Norwell, MA, 2002.

\bibitem{tag} C. Tepedelenlioglu, A. Abdi, and G.B.
Giannakis, ``The Rician $K$ factor: Estimation and performance
analysis,'' {\it IEEE Trans. Wireless Commun.},
vol. 2, pp. 799-810, July 2003.

\bib{valenzuela} R.V. Valenzuela, O. Landron, and D.L. Jacobs,
``Estimating local mean signal strength of indoor multipath
propagation,'' {\it IEEE Trans. Veh. Technol.}, vol. 46, no. 1, Feb.
1997, pp. 203-212.

\end{thebibliography}
\end{document}